\documentclass[10pt]{article}
\usepackage[utf8]{inputenc}
\usepackage{epsfig}
\usepackage{color}
\usepackage[british,english]{babel}
\usepackage{amsthm}
\usepackage{amsmath}
\usepackage{amsfonts}
\usepackage{amssymb}
\usepackage{graphicx}
\newtheorem{corollary}{Corollary}[section]

\newtheorem{lemma}[corollary]{Lemma}
\newtheorem{proposition}[corollary]{Proposition}
\newtheorem{remark}[corollary]{Remark}
\newtheorem{theorem}[corollary]{Theorem}
\newfont{\sBlackboard}{msbm10 scaled 900}

\newcommand{\mylabel}[1]{\label{#1}
            \ifx\undefined\stillediting
            \else \fbox{$#1$}\fi }
\newcommand{\BE}{\begin{equation}}

\newcommand{\EEQ}{\end{equation}}
\newcommand{\rfb}[1]{\mbox{\rm
   (\ref{#1})}\ifx\undefined\stillediting\else:\fbox{$#1$}\fi}

\newfont{\Blackboard}{msbm10 scaled 1200}

\newfont{\roma}{cmr10 scaled 1200}

\def\CC{\rm \hbox{C\kern-.56em\raise.4ex
         \hbox{$\scriptscriptstyle |$}\kern+0.5 em }}

\def \noame{\noalign{\medskip}}
\newcommand{\ep}   {\epsilon}
\usepackage{fancyhdr}

\makeatletter
\def\section{\@startsection {section}{1}{\z@}{-3.5ex plus -1ex minus
    -.2ex}{2.3ex plus .2ex}{\large\bf}}
\makeatother
\def\be{\begin{equation}}
\def\ee{\end{equation}}

\date{ }
\begin{document}
\thispagestyle{empty}
\title{\Large \bf  Navier slip effects in micropolar thin-film flow: a rigorous derivation of Reynolds-type models}\maketitle
\vspace{-2cm}
\begin{center}
Mar\'ia Anguiano\footnote{Departamento de An\'alisis Matem\'atico. Facultad de Matem\'aticas. Universidad de Sevilla. 41012-Sevilla (Spain) anguiano@us.es}, Igor Pa\v zanin\footnote{Department of Mathematics, Faculty of Science, University of Zagreb, Bijeni${\rm \check{c}}$ka 30, 10000 Zagreb (Croatia) pazanin@math.hr} and Francisco J. Su\'arez-Grau\footnote{Departamento de Ecuaciones Diferenciales y An\'alisis Num\'erico. Facultad de Matem\'aticas. Universidad de Sevilla. 41012-Sevilla (Spain) fjsgrau@us.es}
 \end{center}
\ \\
 \renewcommand{\abstractname} {\bf Abstract}
\begin{abstract}
We study the stationary flow of incompressible micropolar fluid in a thin three-dimensional domain under Navier slip boundary condition for the velocity and no-spin condition for microrotation. After rescaling the governing equations, we perform a rigorous asymptotic analysis as the film thickness tends to zero, considering a friction coefficient dependent on the small parameter. According to the scaling of the slip coefficient, we identify three distinct regimes: perfect slip, partial slip, and no-slip. For each regime, we derive the corresponding reduced micropolar system and obtain explicit expressions for the velocity and microrotation fields. This leads to a generalized Reynolds-type equation for the pressure, highlighting the impact of slip effects on the micropolar thin-film flow.
\end{abstract}
\bigskip\noindent

\noindent {\small \bf AMS classification numbers:} 35B27, 35Q35, 76A05.  \\

\ \\
\section {Introduction}\label{S1}
Micropolar fluid models, introduced by Eringen \cite{eringen}, extend the classical Navier-Stokes equations by accounting for the microstructure of the fluid through an additional microrotation field. These models have proven to be effective in describing complex fluids such as lubricants with suspended particles, polymeric fluids, and biological fluids. In lubrication theory, where the flow takes place in a thin domain, micropolar effects can significantly influence the effective pressure and velocity profiles, motivating the derivation of reduced Reynolds-type models (see e.g.~\cite{Luka} and the references therein).\\
\ \\
The asymptotic analysis of micropolar thin-film flows has been extensively performed in the literature. Early rigorous derivations of Reynolds-type equations for micropolar fluids under no-slip boundary conditions were obtained by Bayada and Lukaszewicz \cite{BayadaLuc}, and later extended in various directions, including rough boun\-daries, different scalings of microrotation viscosities, and homogenization effects (see, e.g.~Bayada et al.~\cite{Bayada_Gamouana,Bayada_Gamouana2}, Boukrouche et al.~\cite{Boukrouche3,Boukrouche4}, Pa\v zanin et al.~\cite{Paz1,Paz2}, Anguiano and Suárez-Grau \cite{Anguiano_SG_magneto,Anguiano_SG_Acta}, and the references therein). These works show that the interplay between the domain thickness and the micropolar parameters leads to a variety of effective models, depending on the asymptotic regime.\\
\ \\
In most of the available results on the micropolar thin-film flow, Dirichlet (no-slip) boundary conditions for the velocity field has been imposed. However, from a physical point of view, Navier slip boundary condition for the velocity frequently provides a more realistic description of fluid-boundary interactions, especially when the boundary is smooth or exhibits microscopic roughness. Originally proposed in \cite{Navier}, it is given by
\begin{equation}\label{navierlaw}
[D{\bf u}\, {\bf n}]_{\rm tang}=0,\quad {\bf u}\cdot {\bf n}=0\,,
\end{equation}
where ${\bf n}$ denotes the outward unit normal vector to the boundary, $[D{\bf u}\cdot {\bf n}]_{\rm tang}$ is the tangential component of the vector $D{\bf u} \,{\bf n}$, i.e.~$[D{\bf u}\, {\bf n}]_{\rm tang}:=D{\bf u}\,{\bf n}-[(D{\bf u}\,{\bf n})\cdot {\bf n}]\,{\bf n}$. The boundary condition (\ref{navierlaw}) is often called {\it perfect slip boundary condition} (see \cite{LeRoux}) and it has been rigorously justified by J${\rm \ddot{a}}$ger and Mikeli\'c \cite{Jagger} as homogenized limits of no-slip conditions on rough boundaries. When the boundary is flat, the fluid tends to slip over without friction and there are no boundary layers (see Coron \cite{Coron1, Coron2}), thus Navier boundary condition making more sense than the no-slip condition. In the context of micropolar fluids, the influence of slip boundary conditions has been addressed mainly at a formal or numerical level (see, for instance, Xinhui et al.~\cite{Xinhui}, Bhat and Katagi \cite{Bhat}), while rigorous treatments remain limited. Interesting well-posedness result on non-stationary micropolar flow with Navier slip condition for the velocity can be found in Duarte-Leiva et al.~\cite{Duarte}.\\
\ \\
In the present work, we consider the stationary incompressible linearized micropolar equations
\begin{equation}\label{Micro_intro2_dimension}
\left\{\begin{array}{rl}
\displaystyle
-(\nu+\nu_r)\Delta {\bf  \overline u} +\nabla \overline p=2\nu_r{\rm rot}({\bf\overline  w} )+ {\bf \overline f},\\
\noame
\displaystyle  {\rm div}({\bf \overline u})=0,\\
\noame
\displaystyle -(c_a+c_d)\Delta{\bf \overline w} +4\nu_r{\bf\overline w} =2\nu_r{\rm rot}({\bf \overline w})+ {\bf \overline g},\,
\end{array}\right.
\end{equation}
in a three-dimensional domain given by
\begin{equation}\label{ThinDomain0}
 \overline \Omega^\ep=\left\{\overline x=(\overline x',\overline x_3)\in \mathbb{R}^3\,:\, \overline x'=(\overline x_1, \overline x_2)\in L \omega, \quad 0<\overline x_3< c h\left({\overline x'\over L}\right)
 \right\}\,.
\end{equation}
In view of the application we want to model, the small parameter of the problem is the ratio $\ep =c/L$ between the maximum distance between the surfaces $c$ and $L$, being the characteristic length of the set $\omega\in \mathbb{R}^2$. As usual, ${\bf \overline u}=({\bf \overline u}'(\bar x), \overline u_{3}(\bar x))$, with ${\bf \overline u}'=(\overline u_1, \overline u_2)$, denotes the velocity vector field, $ \overline p= \overline p(\bar x)$ the scalar pressure, and ${\bf  \overline w}=({\bf \overline w}'(\bar x),  \overline w_{3}(\bar x))$, with ${\bf \overline w}'=(\overline w_1, \overline w_2)$, the microrotation field. The given positive constants are the Newtonian viscosity $\nu$ and the microrotation viscosities  $\nu_r, c_a, c_d$ characterizing the isotropic properties of the fluid (see \cite{Luka}). The external forces are given by the functions ${\bf \overline f}=({\bf \overline f}', \overline f_3)$ and ${\bf \overline g}=({\bf \overline g}', \overline g_3)$.\\
\ \\
For the purpose of our analysis, it is convenient to work in dimensionless setting. Thus, introducing the characteristic velocity $V_0$ of
the fluid, and defining
\begin{equation}\label{adimensionalization_change}
\begin{array}{c}
\displaystyle x={\overline x\over L},\quad
{\bf u}={{\bf \overline u}\over V_0},\quad p={L\over  V_0(\nu+\nu_r)}\overline p,\quad {\bf w}={L\over V_0}{\bf \overline w},\quad {\bf f}= {L^2 \over V_0 (\nu+\nu_r)}{\bf \overline f},\quad {\bf g}= {L \over V_0 (\nu+\nu_r)}{\bf \overline g},\\
\noame
\displaystyle
N^2={\nu_r\over \nu+\nu_r},\quad R_M={c_a+c_d\over \nu+\nu_r}{1\over L^2}\,,
\end{array}
\end{equation}
we arrive at
\begin{equation}\label{Micro_intro2}
\left\{\begin{array}{rl}
\displaystyle
-\Delta {\bf u} +\nabla p=2N^2{\rm rot}({\bf w} )+ {\bf  f},\\
\noame
\displaystyle  {\rm div}({\bf u})=0,\\
\noame
\displaystyle -R_M\Delta{\bf w} +4N^2{\bf w} =2N^2{\rm rot}({\bf w})+{\bf g},
\end{array}\right.
\end{equation}
in
\begin{equation}\label{ThinDomain}
 \Omega^\ep=\omega\times (0, \ep h(x')),\quad  x'=(x_1,x_2)\in\omega\subset\mathbb{R}^2,\quad 0<\ep\ll 1\,.
\end{equation}
The dimensionless parameter $N^2<1$ measure the coupling between the velocity and microrotation equations, while $R_M$ is related to characteristic length of the microrotation effects and will be compared with $\ep$ (see \cite{BayadaLuc}).\\
We shall study the system (\ref{Micro_intro2}) endowed with the following type of Navier boundary conditions with no-spin condition on the flat bottom $\Gamma_0=\omega\times \{0\}$:
\begin{equation}\label{navierlaw2}
[D{\bf u}\, {\bf n}]_{\rm tang}=-\lambda[{\bf u}]_{\rm tang},\quad {\bf u}\cdot {\bf n}=0\quad {\bf w}=0\quad \hbox{on }\Gamma_0\,.
\end{equation}
Note that, depending on the value of $\lambda$ (friction coefficient), we have:
\begin{itemize}
\item Perfect slip when $\lambda=0$, i.e. condition (\ref{navierlaw}),
\begin{equation}\label{navierlaw2perfect}
[D{\bf u}\, {\bf n}]_{\rm tang}=0,\quad {\bf u}\cdot {\bf n}=0\quad {\bf w}=0\quad \hbox{on }\Gamma_0,
\end{equation}
\item Partial slip when $\lambda\in (0,+\infty)$, i.e.
\begin{equation}\label{navierlaw2partial}
[D{\bf u}\, {\bf n}]_{\rm tang}=-\lambda[{\bf u}]_{\rm tang},\quad {\bf u}\cdot {\bf n}=0\quad {\bf w}=0\quad \hbox{on }\Gamma_0,
\end{equation}
\item No-slip when  $\lambda=+\infty$,
\begin{equation}\label{navierlaw2}
[{\bf u}]_{\rm tang}=0,\quad {\bf u}\cdot {\bf n}=0\quad {\bf w}=0\quad \hbox{on }\Gamma_0.
\end{equation}
\end{itemize}
Assuming no-slip and no-spin conditions on the upper boundary $\Gamma_1^\ep=\omega\times \{\ep h(x')\}$, the goal of this paper is to identify the effects of the Navier slip boundary conditions on the effective behavior of the fluid flow. In addition, the friction coefficient in the Navier condition is allowed to depend on the domain's thickness $\epsilon$ (see (\ref{lambdaep}) below), which enables us to capture different slip regimes within a unified framework.\\
\ \\
By means of compactness arguments and careful asymptotic analysis, we identify three distinct limit regimes depending on the scaling of the friction coefficient: perfect slip, partial slip, and no-slip. For each regime, we rigorously derive the corresponding reduced effective system, which consists of a coupled one-dimensional micropolar problem in the transverse variable and a two-dimensional Reynolds-type equation for the pressure. The limit velocity and microrotation profiles are explicitly computed, allowing us to determine how slip effects modify the effective permeability and the structure of the generalized Reynolds equation. From a methodological point of view, the presented analysis follows the asymptotic framework developed for thin-domain fluid problems and micropolar thin-film models. A priori estimates for the velocity, microrotation, and pressure fields are derived using energy methods adapted to Navier slip boundary conditions, in the spirit of Bayada and Lukaszewicz \cite{BayadaLuc} and the subsequent works \cite{Bayada_Gamouana,Bayada_Gamouana2,SG1}. The pressure is handled through a decomposition technique in thin domains, relying on compactness and duality arguments as in \cite{CLS,grau1}.\\
\ \\
The paper is organized as follows. In Section 2, after formally describing the problem under consideration, we prove the variational formulation and the corresponding  existence and uniqueness result. Section 3 is devoted to the rescaling of the problem into a fixed reference domain, which is essential for performing the asymptotic analysis. In Section 4, we derive uniform a priori estimates for the unknowns and establish the main compactness and convergence results as the thickness parameter tends to zero. The limit problem is rigorously identified in Section 5, where we characterize the effective reduced system and determine the boundary conditions in the limit depending on the scaling of the slip coefficient. Explicit expressions for the limit velocity and microrotation profiles are obtained in Section 6, allowing a detailed analysis of the flow structure. Finally, in Section 7, we derive the generalized Reynolds equation for the pressure and discuss the influence of Navier slip effects on the resulting effective model.

\section{Problem setting and the well-posedness result}\label{sec:setting}
\subsection{Description of the domain and some notation}  Let $\omega\in \mathbb{R}^2$ be a connected set with a Lipschitz boundary. As indicated in the Introduction, the considered thin domain is defined by
\begin{equation}\label{Omegaep}
\Omega^\epsilon=\{x=(x',x_3)\in\mathbb{R}^2\times \mathbb{R}\,:\, x'\in \omega,\ 0<x_3< h_\epsilon(x')\}\,.
\end{equation}
Here the function $h_\ep(x')= \epsilon h\left(x'\right)$ represents the real distance between the two surfaces, where the function $h$ is $W^{1,\infty}$ and positive, defined for $x'$. Furthermore, we assume that there exist $h_{\rm min}$ and $h_{\rm max}$ such that $0<h_{\rm min}\leq h(x')<h_{\rm max}$ for all $(x',0)\in \omega$.\\
\ \\
We denote the boundaries of $\Omega^\epsilon$ as follows
$$\begin{array}{c}
\displaystyle \Gamma_0 =\omega\times\{0\},\quad \Gamma_1^\epsilon=\left\{(x',x_3)\in\mathbb{R}^2\times\mathbb{R}\,:\, x'\in \omega,\ x_3=  h_\epsilon(x')\right\},\quad
\displaystyle\Gamma_{\rm lat}^\epsilon=\partial\Omega_\ep\setminus (\Gamma_0\cup\Gamma_1^\ep).
\end{array}$$

\noindent Employing a dilatation in the vertical variable, namely. $z_3=x_3/\ep$, w introduce the following rescaled sets ($z'=x'$):
 \begin{equation}\label{domains_tilde}\begin{array}{c}
 \displaystyle   \Omega =\{(z',z_3)\in\mathbb{R}^2\times \mathbb{R}\,:\, z'\in \omega,\ 0<z_3< h(z')\},\\
 \noame
 \Gamma_1 =\{(z',z_3)\in\mathbb{R}^2\times \mathbb{R}\,:\, z'\in \omega,\ z_3=h(x')\}\quad   \Gamma_{\rm lat}=\partial\Omega \setminus (\Gamma_0\cup  \Gamma_1).
  \end{array}
  \end{equation}

\noindent From now on, we denote by $C$ a generic constant, independent of $\ep$. Moreover, $O_\ep$ denotes a  generic quantity, which can change from line to line, meant to tend to zero when $\ep\to 0$. Finally, we employ the usual notation for the partial differential operators:
 $$\begin{array}{c}
 \displaystyle \Delta{\bf \varphi}=\Delta_{x'}{\bf \varphi} +\partial_{x_3}^2 {\bf \varphi},\quad {\rm div}({\bf \varphi} )={\rm div}_{x'}(\varphi')+\partial_{x_3}\varphi_3,\\
 \\
  \displaystyle  {\rm rot}({\bf \varphi} )= (\partial_{x_2}\varphi_3-\partial_{x_3}\varphi_2,
  -\partial_{x_1}\varphi_3+\partial_{x_3}\varphi_1,\partial_{x_1}\varphi_2-\partial_{x_2}\varphi_1)^t,
 \end{array}$$
where ${\bf \varphi} =(\varphi',\varphi_3)$ with $\varphi'=(\varphi_1, \varphi_2)$, is a vector  function defined in $\Omega^\epsilon$. \\

\subsection{Micropolar equations with Navier slip boundary conditions} Assuming small Reynolds number and neglecting the inertia effects, we study (dimensionless) linearized micro\-polar equations in stationary regime:
\begin{equation}\label{system_1}
\left\{\begin{array}{rl}
\displaystyle -\Delta {\bf u}_\ep+\nabla p_\ep=2N^2{\rm rot}( {\bf w}_\ep)+{\bf f}_\ep & \hbox{in}\ \Omega^\epsilon,\\
\\
{\rm div}( {\bf u}^\epsilon)=0& \hbox{in}\ \Omega^\epsilon,\\
\\
\displaystyle -R_M\Delta  {\bf w}_\ep+4N^2{\bf w}_\ep=2N^2{\rm rot}({\bf u}_\ep)+{\bf g}_\ep& \hbox{in}\ \Omega^\epsilon\,.
\end{array}\right.
\end{equation}
In the above system, the unknown functions are
$${\bf u}_\ep=({\bf u}'_\ep(x), u_{3,\ep}(x)),\quad {\bf w}_\ep=({\bf w}'_\ep(x), w_{3,\ep}(x)),\quad p_\ep=p_\ep(x),$$
representing the velocity, the microrotation and the
pressure of the fluid respectively at a point $x\in \Omega^\ep$. Note that we add a subscript $\ep$ to stress the dependence of the solution on the small parameter.\\
\ \\
Assuming that the mean value of the pressure is zero, we prescribe the following boundary conditions on the top and lateral boundaries
\begin{equation}\label{BCBot}
\begin{array}{ll}
\displaystyle {\bf u}_\epsilon=0,\quad  {\bf w}_\epsilon=0&\hbox{on}\ \Gamma_1^\ep\cup {\Gamma_{\rm lat}^\ep},
\end{array}
\end{equation}
and the Navier slip boundary conditions with no-spin condition on the bottom of the domain:
\begin{equation}\label{BCTopBot1}
[D{\bf u}_\ep\, {\bf n}]_{\rm tang}=-\lambda_\ep [{\bf u}_\ep]_{\rm tang},\quad {\bf u}_\ep\cdot {\bf n}=0,\quad {\bf w}_\ep=0\quad \hbox{on }\Gamma_0\,.
\end{equation}
In our analysis, we allow the friction coefficient $\lambda_\ep$ to depend on $\ep$, i.e.
\begin{equation}\label{lambdaep}\lambda_\ep=\lambda\,\ep^\gamma,\quad \lambda>0,\quad \gamma\in \mathbb{R}.
\end{equation}
\noindent It should be observed that the order of magnitude of the dimensionless parameter $R_M$ plays an important role with regards to the asymptotic behavior of the solution (see e.g. \cite{Bayada_Gamouana, BayadaLuc, Pazanintaml}). We consider the most interesting one, which  leads to a strong
coupling at main order, assuming
\begin{equation}\label{RMnumber}R_M=\epsilon^2\,R_c,\quad R_c>0.
\end{equation}
According to \cite{grau1, SG1}, we also consider the following assumptions on the forces
$${\bf f}_\ep(x)=({\bf f}'(x'), 0)^t,\quad {\bf g}_\ep(x)=\ep({\bf g}'(x'), 0)^t,$$
which means that, due to the thickness of the domain, the vertical components of the forces can be neglected and the forces can be considered independent of the vertical variable.\\
\ \\
The main goal of the present paper is to study the asymptotic behavior of the solution to the problem (\ref{system_1})--(\ref{RMnumber}) via rigorous asymptotic analysis with respect to the small parameter $\ep$. However, before that, we need to establish the corresponding existence and uniqueness result. In view of that, we introduce the notion of weak solution to system (\ref{system_1})--(\ref{BCTopBot1}). Taking into account the boundary conditions, it is natural to consider the following functional spaces
$$\begin{array}{c}
\displaystyle V^\ep=\{\varphi\in H^1(\Omega^\ep)^3\,:\, \varphi=0\quad \hbox{in }\Gamma_1^\ep\cup \Gamma_{\rm lat}^\ep,\quad \varphi\cdot {\bf n}=\varphi_3=0\quad\hbox{ on }\Gamma_0\},\quad V_0^\ep=\{\varphi\in V^\ep\,:\, {\rm div}(\varphi)=0\},
\end{array}$$
equipped with the norm of $\|D\varphi\|_{L^2(\Omega^\ep)^{3\times 3}}$, and the space
$$L^2_0(\Omega^\ep)=\left\{q\in L^2(\Omega^\ep)\,:\, \int_{\Omega_\ep}q\,dx=0\right\},$$
equipped with the norm of $L^2(\Omega^\ep)$. Since $\Gamma_0$ is flat, ${\bf n}=(0,0,-1)^t$ and  $[{\bf u}_\ep]_{\rm tang}={\bf u}'_\ep$ on $\Gamma_0$, it holds:
\begin{equation}\label{Parts1}
\begin{array}{rl}
\displaystyle \int_{\Omega^\ep}{\rm rot}({\bf u}_\ep)\cdot \psi\,dx=&\displaystyle\int_{\Omega^\ep}{\rm rot}(\psi)\cdot {\bf u}_\ep\,dx-\int_{\Gamma_0}({\bf u}_\ep \times {\bf n})\cdot \psi\,d\sigma(x')\\
\noame
=&\displaystyle\int_{\Omega^\ep}{\rm rot}(\psi)\cdot {\bf u}_\ep\,dx\quad\quad\forall \psi\in H^1_0(\Omega^\ep)^3,
\end{array}
\end{equation}
\begin{equation}\label{Parts2}
\begin{array}{rl}
\displaystyle \int_{\Omega^\ep}{\rm rot}({\bf w}_\ep)\cdot \varphi\,dx=&\displaystyle\int_{\Omega^\ep}{\rm rot}(\varphi)\cdot {\bf w}_\ep\,dx-\int_{\Gamma_0}({\bf w}_\ep \times {\bf n})\cdot \varphi\,d\sigma(x')\\
\noame
=&\displaystyle\int_{\Omega^\ep}{\rm rot}(\varphi)\cdot {\bf w}_\ep\,dx\quad\quad\forall \varphi\in V^\ep(\Omega^\ep)^3,
\end{array}
\end{equation}
\begin{equation}\label{Parts2}
\begin{array}{rl}
\displaystyle -\int_{\Omega^\ep}\Delta {\bf u}_\ep\cdot \varphi\,dx=&\displaystyle\int_{\Omega^\ep}D{\bf u}_\ep: D \varphi\,dx-\int_{\Gamma_0}[D{\bf u}_\ep\,{\bf n}]_{\rm tang}\cdot  \varphi'\,d\sigma(x')\\
\noame
=&\displaystyle
\displaystyle\int_{\Omega^\ep}D{\bf u}_\ep: D \varphi\,dx+\lambda\ep^\gamma\int_{\Gamma_0}{\bf u}_\ep'\cdot  \varphi'\,d\sigma(x')\quad\quad\forall \varphi\in V^\ep(\Omega^\ep).
\end{array}
\end{equation}
According to properties (\ref{Parts1})-(\ref{Parts2}), a weak solution to (\ref{system_1})--(\ref{BCTopBot1}) is defined as follows:

Find a triplet $({\bf u}_\ep, {\bf w}_\ep, p_\ep)\in V^\ep_0\times H^1_0(\Omega^\ep)\times L^2_0(\Omega^\ep)$ such that
\begin{eqnarray}
&&\displaystyle \int_{\Omega_\ep}D{\bf u}_\ep : D\varphi\,dx-\int_{\Omega_\ep}p_\ep\,{\rm div}(\varphi)\,dx
-2N^2\int_{\Omega_\ep} {\rm rot}(\varphi)\cdot{\bf w}_\ep\,dx\nonumber\\
\noame
&&\displaystyle   +\lambda\ep^\gamma\int_{\Gamma_0}{\bf u}'_\ep\cdot \varphi'\,d\sigma(x')=\int_{\Omega_\ep}{\bf f}'\cdot \varphi'\,dx,\label{Form_Var_vel}
\\
\noame
&&\displaystyle \ep^2 R_c\int_{\Omega_\ep}D{\bf w}_\ep:D\psi\,dx+4N^2\int_{\Omega_\ep}{\bf w}_\ep\cdot\psi\,dx-2N^2\int_{\Omega_\ep}{\bf u}_\ep\cdot {\rm rot}(\psi)\,dx= \ep \int_{\Omega_\ep}{\bf g}' \cdot \psi'\,dx, \label{Form_Var_micro}
\end{eqnarray}
\indent for every $\varphi\in V^\ep$ and $\psi\in H^1_0(\Omega^\ep)^3$.
\begin{theorem}\label{thm:existence}
Under the previous assumptions, for every $\ep>0$,  there exists a unique weak solution $({\bf u}_\epsilon, {\bf w}_\epsilon, p_\epsilon)\in V^\ep_0\times H^1_0(\Omega^\ep)^3\times L^2_0(\Omega^\ep)$ of (\ref{Form_Var_vel})--(\ref{Form_Var_micro}).
\end{theorem}
\begin{proof} We divide the proof in two steps.\\

{\it Step 1.  Mixed variational formulation} It easy to see that $({\bf u}_\ep, {\bf w}_\ep, p_\ep)$ is a weak solution to system (\ref{system_1})--(\ref{BCTopBot1}) if and only if $({\bf u}_\ep, {\bf w}_\ep, p_\ep)$ satisfies for any $(\varphi,\psi)\in V^\ep\times H^1_0(\Omega^\ep)^3$:
\begin{equation}\label{formvarsum}
\begin{array}{l}\displaystyle\int_{\Omega_\ep}D{\bf u}_\ep : D\varphi\,dx-\int_{\Omega_\ep}p_\ep\,{\rm div}(\varphi)\,dx
-2N^2\int_{\Omega_\ep} {\rm rot}(\varphi)\cdot{\bf w}_\ep\,dx   +\lambda\ep^\gamma\int_{\Gamma_0}{\bf u}'_\ep\cdot \varphi'\,d\sigma(x')\\
\noame
\displaystyle+\ep^2 R_c\int_{\Omega_\ep}D{\bf w}_\ep:D\psi\,dx+4N^2\int_{\Omega_\ep}{\bf w}_\ep\cdot\psi\,dx-2N^2\int_{\Omega_\ep}{\bf u}_\ep\cdot {\rm rot}(\psi)\,dx\\
\noame
\displaystyle
=\int_{\Omega_\ep}{\bf f}'\cdot \varphi'\,dx+\ep \int_{\Omega_\ep}{\bf g}'\cdot \psi'\,dx.
\end{array}
\end{equation}
Equation (\ref{formvarsum}) justifies the introduction of the bilinear form $\mathcal{A}_\ep: V^\ep\times H_0^1(\Omega^\ep)^3\to \mathbb{R}$ and $\mathcal{B}_\ep: (V^\ep\times H_0^1(\Omega^\ep)^3)\times L^2_0(\Omega^\ep) \to \mathbb{R}$ respectively defined by
\begin{equation}\label{A_bilinear_1}
\begin{array}{rl}
\mathcal{A}_\ep(({\bf  u}, {\bf   w}), (  \varphi, \psi))=& \displaystyle  \int_{\Omega_\ep}D{\bf u} : D\varphi\,dx-\int_{\Omega_\ep}p_\ep\,{\rm div}(\varphi)\,dx
-2N^2\int_{\Omega_\ep} {\rm rot}(\varphi)\cdot{\bf w}\,dx  +\lambda\ep^\gamma \int_{\Gamma_0}{\bf u}'\cdot \varphi'\,d\sigma(x'),\\
\noame
&\displaystyle+\ep^2 R_c\int_{\Omega_\ep}D{\bf w}:D\psi\,dx+4N^2\int_{\Omega_\ep}{\bf w}\cdot\psi\,dx-2N^2\int_{\Omega_\ep}{\bf u}\cdot {\rm rot}(\psi)\,dx,
\end{array}
\end{equation}
\begin{equation}\label{B_linear}
\mathcal{B}^\ep(( {\bf u}, {\bf w}), p_\ep)=- \int_{ \Omega^\ep}  p_\ep\,{\rm div}({\bf v})\,dx,
\end{equation}
and $\mathcal{L}^\ep: V^\ep\times H_0^1(\Omega)^3\to \mathbb{R}$ given by
\begin{equation}\label{L_linear}
\displaystyle\mathcal{L}^\ep({\bf u}, {\bf w})= \int_{\Omega_\ep}{\bf f}'\cdot {\bf v}'\,dx+\ep\int_{\Omega_\ep}{\bf g}'\cdot {\bf w}'\,dx.
\end{equation}
Therefore, $({\bf u}_\ep, {\bf w}_\ep, p_\ep)$ is a weak solution of system (\ref{system_1})--(\ref{BCTopBot1})  if and only if $({\bf u}_\ep, {\bf w}_\ep, p_\ep)\in V^\ep\times H_0^1(\Omega^\ep)^3\times L^2_0(\Omega^\ep)$ and satisfies the following mixed formulation
\begin{eqnarray}
\mathcal{A}^\ep(({\bf u}_\ep, {\bf  w}_\ep), (  \varphi, \psi))+\mathcal{B}^\ep((  \varphi,\ \psi),  p_\ep)=\mathcal{L}( \varphi, \psi)&\forall\,(  \varphi,  \psi)\in   V^\ep\times  H^1_0(\Omega^\ep)^3,\label{Form_Var_A}\\
\noame
\mathcal{B}^\ep(({\bf   u}_\ep, {\bf   w}_\ep),  q)=0&\forall\,   q\in L^2_0(\Omega^\ep).
\end{eqnarray}

{\it Step 2.  Existence and uniqueness of solution of the mixed variational formulation.} From the Cauchy-Schwarz inequality and the Poincar\'e inequality in a thin domain of thickness $\ep$ (see Lemma \ref{Poincare_lemma} below), i.e.
$$\|\varphi\|_{L^2(\Omega^\ep)^3}\leq C\ep \|D\varphi\|_{L^2(\Omega^\ep)^{3\times 3}},$$
for any $\varphi\in V^\ep$ or  $\varphi\in H^1_0(\Omega^\ep)^3$, it is easy to see that $\mathcal{A}_\ep$, $\mathcal{B}_\ep$ and $\mathcal{L}_\ep$ are continuous on their respective domains of definitions, for any fixed value of parameter $\ep$. Therefore, taking into account that by definition of $V^\ep_0$ it holds
$$V^\ep_0\times H_0^1(\Omega^\ep)^3=\{(\varphi,\psi)\in V^\ep_0\times H_0^1(\Omega^\ep)^3\ :\ \mathcal{B}_\ep((\varphi, \psi), q)=0\ \hbox{for any }q\in L^2_0(\Omega^\ep)\},$$
and denoting the norm
$$\|(\varphi,\psi)\|_{V^\ep\times H^1_0}=\left(\|D\varphi\|_{L^2(\Omega^\ep)^{3\times 3}}+\|D\psi\|_{L^2(\Omega^\ep)^{3\times 3}}\right)^{1/2},$$
the existence and uniqueness of the solution $({\bf u}_\ep, {\bf w}_\ep, p_\ep)$ to the mixed variational formulation result from the following properties (see \cite[paragraph 4.1 p. 57]{Girault} for more details):
\begin{itemize}
\item[(i)] Coerciveness of $\mathcal{A}_\ep$: there exists $\eta=\eta(\ep)>0$ such that
$$ \mathcal{A}_\ep((\varphi, \psi), (\varphi, \psi))\geq \eta \|(\varphi,\psi)\|_{V^\ep\times H^1_0}^2,\quad \forall(\varphi, \psi)\in V^\ep_0\times H^1_0(\Omega^\ep)^3,$$
\item[(ii)] Inf-sup condition: there exists $c=c(\ep)$ such that
$$\inf_{q\in L^2_0(\Omega^\ep)}\sup_{(\varphi, \psi)\in V^\ep\times H^1_0(\Omega^\ep)^3}{\mathcal{B}_\ep((\varphi,\psi),p)\over\|(\varphi,\psi)\|_{V_\ep\times H_0^1}\|q\|_{L^2_0(\Omega^\ep)}}\geq c.$$
\end{itemize}
To proof the coerciveness condition $(i)$, we use the Cauchy-Schwarz inequality and so
$$\begin{array}{rl}
\mathcal{A}_\ep((  \varphi, \psi), (  \varphi, \psi))=& \displaystyle  \int_{\Omega_\ep}|D\varphi|^2\,dx
+\lambda\ep^\gamma \int_{\Gamma_0}|\varphi'|^2\,d\sigma(x')\\
\noame
&\displaystyle+\ep^2 R_c\int_{\Omega_\ep}|D\psi|^2\,dx+4N^2\int_{\Omega_\ep}|\psi|^2\,dx-4N^2\int_{\Omega_\ep}{\rm rot}(\varphi)\cdot  \psi\,dx,\\
\noame
\geq & \displaystyle \|D\varphi\|_{L^2(\Omega^\ep)^{3\times 3}}^2+\ep^2 R_c\|D\psi\|_{L^2(\Omega^\ep)^{3\times 3}}^2+4N^2\|\psi\|_{L^2(\Omega^\ep)^3}^2\,dx-4N^2\|D\varphi\|_{L^2(\Omega^\ep)^{3\times 3}}\|\psi\|_{L^2(\Omega^\ep)^3}
\end{array}
$$
By using the Young inequality, we deduce
$$4N^2\|D\varphi\|_{L^2(\Omega^\ep)^{3\times 3}}\|\psi\|_{L^2(\Omega^\ep)^3}
\leq N^2\|D\varphi\|_{L^2(\Omega^\ep)^{3\times 3}}^2+4N^2\|\psi\|^2_{L^2(\Omega^\ep)^3},$$
and so, we get
$$\begin{array}{rl}
\mathcal{A}_\ep((  \varphi, \psi), (  \varphi, \psi))\geq & \displaystyle
(1-N^2)\|D\varphi\|_{L^2(\Omega^\ep)^{3\times 3}}^2+\ep^2 R_c\|D\psi\|_{L^2(\Omega^\ep)^{3\times 3}}^2,
\end{array}
$$
which implies the coerciveness of $\mathcal{A}_\ep$ with $\eta=\min\{1-N^2, \ep^2 R_c\}$.\\

To prove $(ii)$, let $(\varphi,\psi)$ belong to $V^\ep\times H_0^1(\Omega^\ep)^3$, and $q\in L^2_0(\Omega^\ep)$. Then, we obtain
$$H_0^1(\Omega^\ep)^3\times \{0\}\subset V^\ep\times H_0^1(\Omega^\ep)^3,$$
so that
\begin{equation}\label{infsup}\sup_{(\varphi,\psi)\in V^\ep\times H_0^1(\Omega^\ep)^3}{\int_{\Omega^\ep}{\rm div}(\varphi)q\,dx\over (\|\varphi\|_{V^\ep}+\|\psi\|_{H_0^1})^{1/2}}\geq \sup_{(\varphi,0)\in   H_0^1(\Omega^\ep)^3\times \{0\}}{\int_{\Omega^\ep}{\rm div}(\varphi)q\,dx\over \|D\varphi\|_{L^2(\Omega^\ep)^3}}=\sup_{\varphi \in   H_0^1(\Omega^\ep)^3}{\int_{\Omega^\ep}{\rm div}(\varphi)q\,dx\over \|D\varphi\|_{L^2(\Omega^\ep)^3}}.
\end{equation}
According to \cite[Lemma 3]{Duvnjak}, for a given $q\in L^2_0(\Omega^\ep)$, there is $\varphi_0=\varphi(q)$ such that
$${\rm div}(\varphi_0)=q\quad \hbox{and}\quad \|D\varphi_0\|_{L^2(\Omega^\ep)^{3\times 3}}\leq {C\over \ep}\|q\|_{L^2(\Omega^\ep)}.$$  Choosing $\varphi=\varphi_0$ in (\ref{infsup}), we deduce
$$\sup_{(\varphi,\psi)\in V^\ep\times H_0^1(\Omega^\ep)^3}{\int_{\Omega^\ep}{\rm div}(\varphi)q\,dx\over (\|\varphi\|_{V^\ep}+\|\psi\|_{H_0^1})^{1/2}}\geq  {\int_{\Omega^\ep}{\rm div}(\varphi_0)q\,dx\over \|D\varphi_0\|_{L^2(\Omega^\ep)^3}}={\|q\|^2_{L^2(\Omega^\ep)}\over \|D\varphi_0\|_{L^2(\Omega^\ep)^{3\times 3}}}\geq {\ep\over C}\|q\|_{L^2(\Omega^\ep)},$$
which implies condition $(ii)$ with $c(\ep)={\ep\over C}$. This ends the proof.
\end{proof}

\section{Rescaled problem}\label{sec:rescaling}
In order to deduce the asymptotic behavior of the solution (${\bf u}_\ep, {\bf w}_\ep$, $p_\ep$) when $\ep\rightarrow 0$, we rescale the governing problem by applying the dilatation in the variable $x_3$ given by
\begin{equation}\label{dilatacion}
\quad z_3={x_3\over \ep}\,.
\end{equation}
Denoting $z'=x'$, consequently we obtain the functions defined in  $\Omega$ and boundaries $\Gamma_1$ and  $\Gamma_{\rm lat}$ (see (\ref{domains_tilde})) satisfying:\\
\begin{equation}\label{system_1_dil}
\left\{\begin{array}{rl}
\displaystyle -\Delta_\ep {\bf \widetilde u}_\ep +\nabla_\ep \widetilde p_\ep=2N^2{\rm rot}_\ep( {\bf \widetilde w}_\ep)  +{\bf f}'& \hbox{in}\  \Omega ,\\
\\
{\rm div}_\ep( {\bf \widetilde v}^\epsilon)=0& \hbox{in}\  \Omega,\\
\\
\displaystyle -\ep^2 R_c\Delta_\ep   {\bf \widetilde w}_\ep+4N^2 {\bf \widetilde w}_\ep=2N^2{\rm rot}_\ep({\bf \widetilde u}_\ep)+\ep {\bf g}'& \hbox{in}\  \Omega\,,
\end{array}\right.\end{equation}
with
\begin{equation}\label{BCBottilde}
\begin{array}{l}
\displaystyle \widetilde {\bf  u}_\epsilon=0,\quad \widetilde {\bf w}_\epsilon=0\quad \hbox{on}\  \Gamma_1\cup {\Gamma_{\rm lat}},
\\
\noame
\displaystyle
[D_\ep {\bf \widetilde u}_\ep\, {\bf n}]_{\rm tang}=-\lambda\ep^\gamma {\bf \widetilde u}_\ep',\quad {\bf \widetilde  u}_{\ep,3}=0\quad {\bf \widetilde  w}_\ep=0\quad \hbox{on }\Gamma_0\,.
\end{array}
\end{equation}
Here ${\bf \widetilde u}_\epsilon(z',z_3)={\bf u}^\epsilon(z',\epsilon z_3)$, $ \widetilde p_\epsilon(z',z_3)=p_\epsilon(z',\epsilon z_3)$, ${\bf \widetilde w}_\epsilon(xz,z_3)={\bf w}_\epsilon(z',\epsilon z_3)$  for a.e. $z\in  \Omega$. The rescaled operators are given by
\begin{equation}\label{def_operator_tilde}\begin{array}{c}
 \displaystyle \Delta_{\epsilon}{\bf \widetilde\varphi}=\Delta_{z'}{\bf \widetilde\varphi}+{1\over \ep^2}\partial_{z_3}^2{\bf \widetilde\varphi},\quad {\rm div}_{\epsilon}({\bf \widetilde \varphi})={\rm div}_{z'}(\widetilde \varphi')+{1\over \ep}\partial_{z_3}\widetilde \varphi_3,
 \\
 \noame
\displaystyle  {\rm rot}_{\epsilon}({\bf\widetilde  \varphi})=\left( {\rm rot}_{z'}(\widetilde\varphi_3)+{1\over \ep}{\rm rot}_{z_3}(\widetilde\varphi'), {\rm Rot}_{z'}(\widetilde\varphi')\right)^t,
 \end{array}
 \end{equation}
for a vector function $\widetilde\varphi=(\widetilde\varphi',\widetilde\varphi_3)$ defined in $\Omega$. Hereinafter, denoting by
 \begin{equation}\label{vec_perp}(\widetilde\varphi')^\perp=(-\widetilde\varphi_2,\widetilde\varphi_1)^t,
 \end{equation}
 we use the following notation
 \begin{equation}\label{rotationals}
 {\rm rot}_{z'}(\widetilde\varphi_3)=(\partial_{z_2}\widetilde\varphi_3,-\partial_{z_1}\widetilde\varphi_3)^t,\quad {\rm rot}_{z_3}(\widetilde\varphi')=\partial_{z_3}(\widetilde\varphi')^\perp,\quad {\rm Rot}_{z'}(\widetilde\varphi')=\partial_{z_1}\widetilde\varphi_2-\partial_{z_2}\widetilde\varphi_1.
 \end{equation}
Now, the rescaled variational formulation in $\Omega$ reads

Find $({\bf \widetilde v}_\ep, {\bf \widetilde w}_\ep, \widetilde p_\ep)\in \widetilde V_0\times H^1_0(\Omega)^3\times L^2_0(\Omega)$ such that
\begin{equation}\label{Form_Var_vel_tilde}\begin{array}{l}
\displaystyle  \int_{\Omega}D_\ep{\bf \widetilde u}_\ep : D_\ep\widetilde \varphi\,dz- \int_{\Omega}\widetilde p_\ep\,{\rm div}_\ep(\widetilde \varphi)\,dz
-2N^2 \int_{\Omega} {\rm rot}_\ep(\widetilde \varphi)\cdot{\bf \widetilde w}_\ep\,dz \\
\noame
\displaystyle   +\lambda\ep^{\gamma-1}\int_{\Gamma_0}{\bf \widetilde u}'_\ep\cdot \widetilde \varphi'\,d\sigma(z')= \int_{\Omega}{\bf f}'\cdot \widetilde \varphi'\,dz,
\\
\noame
\displaystyle \ep^2 R_c\int_{\Omega}D_\ep{\bf \widetilde w}_\ep:D_\ep\widetilde \psi\,dz+4N^2 \int_{\Omega}{\bf\widetilde  w}_\ep\cdot\widetilde \psi\,dz-2N^2\int_{\Omega}{\bf \widetilde u}_\ep\cdot {\rm rot}_\ep(\widetilde \psi)\,dz= \ep \int_{\Omega}{\bf g}' \cdot \widetilde \psi'\,dz,
\end{array}
\end{equation}
for every $(\widetilde \varphi,\widetilde\psi)\in \widetilde V\times H^1_0(\Omega)^3$ obtained from $(\varphi, \psi)$ via (\ref{dilatacion}), where
$$\begin{array}{c}
\displaystyle \widetilde V=\{\widetilde \varphi\in H^1(\Omega)^2\,:\, \widetilde \varphi=0\quad \hbox{in }  \Gamma_1\cup  \Gamma_{\rm lat},\quad \widetilde \varphi_3=0\quad\hbox{ on }\Gamma_0\},\quad \widetilde V_0=\{\widetilde \varphi\in \widetilde V\,:\, {\rm div}_\ep(\widetilde \varphi)=0\}.
\end{array}
$$

\section{{\it A priori} estimates and convergences} \label{sec:estimates}
In this section, we derive the {\it a priori} estimates and obtain some compactness results. First, we recall some auxiliary results that will be useful in the sequel. Then we prove the {\it a priori} estimates for velocity and microrotation in the $\ep$-independent domain $\Omega$. By using these estimates, we derive the corresponding {\it a priori} estimates for pressure in $\Omega$ by means a decomposition result for the pressure. Finally, we employ the derived {\it a priori} estimates, to obtain the convergence results.

\subsection{Auxiliary results}
To derive the desired estimates, let us recall some known technical results (see e.g.~\cite{grau1}).
\begin{lemma}\label{Poincare_lemma} For any $\varphi\in V^\ep$ or $H_0^1(\Omega^\ep)^3$, there holds the following inequality
\begin{equation}\label{Poincare}
\|\varphi\|_{L^2(\Omega^\epsilon)^3}\leq C\epsilon\|D \varphi\|_{L^2(\Omega^\epsilon)^{3\times 3}}.
\end{equation}
Moreover, from the change of variables (\ref{dilatacion}),   there holds for all $\widetilde\varphi\in \widetilde V$  or $H_0^1(\Omega)^3$  the following rescaled estimate
\begin{equation}\label{Poincare2}
\|\widetilde \varphi\|_{L^2( \Omega)^3}\leq C\epsilon\|D_{\ep} \widetilde \varphi\|_{L^2(  \Omega)^{3\times 3}}.
\end{equation}
\end{lemma}
\begin{lemma}
For all $\varphi \in V^\ep$ or $H_0^1(\Omega^\ep)^3$, it holds
\begin{equation}\label{Gaffney}
 \|D\varphi\|_{L^2(\Omega^\ep)^{3\times 3}}^2=\|{\rm rot}(\varphi)\|^2_{L^2(\Omega^\ep)^{3}}+\|{\rm div}(\varphi)\|^2_{L^2(\Omega^\ep)}, \quad \|{\rm rot}(\varphi)\|^2_{L^2(\Omega^\ep)^{3}}\leq \|D\varphi\|_{L^2(\Omega^\ep)^{3\times 3}}^2,
\end{equation}
and, if moreover, ${\rm div}(\varphi)=0$ in $\Omega^\ep$, then
\begin{equation}\label{Gaffney_div0}
\|{\rm rot}(\varphi)\|_{L^2(\Omega^\ep)^{3}}^2=\|D\varphi\|_{L^2(\Omega^\ep)^{3\times 3}}^2.
\end{equation}
Moreover, from the change of variables (\ref{dilatacion}), we have for all $\widetilde\varphi\in \widetilde V$ or $H_0^1(\Omega)^3$ that
\begin{equation}\label{Gaffney_tilde}\|D_\ep \widetilde \varphi\|^2_{L^2( \Omega)^{3\times 3}}=\|{\rm rot}_\ep(\widetilde \varphi)\|^2_{L^2(\Omega)^{3}}+\|{\rm div}_\ep(\widetilde \varphi)\|^2_{L^2( \Omega)},\quad \|{\rm rot}_\ep(\widetilde \varphi)\|^2_{L^2(  \Omega)^{3}}\leq \|D_\ep \widetilde \varphi\|^2_{L^2(  \Omega)^{3\times 3}},
\end{equation}
and if moreover, ${\rm div}_\ep(\widetilde \varphi)=0$ in $ \Omega$, then
\begin{equation}\label{Gaffney_div0_tilde}
\|{\rm rot}_\ep\widetilde \varphi\|_{L^2(\Omega)^{3}}^2=\|D_\ep\widetilde \varphi\|_{L^2(\Omega)^{3\times 3}}^2.
\end{equation}
\end{lemma}
\begin{proof} For all $\varphi\in V^\ep$, it holds (see for instance \cite{Boyer} formula (IV.23))
$$\int_{\Omega^\ep}(|{\rm div}(\varphi)|^2+|{\rm rot}(\varphi)|^2)\,dx=\int_{\Omega^\ep}|D\varphi|^2\,dx+\int_{\Gamma_0}((\varphi\cdot\nabla){\bf n})\cdot \varphi\,d\sigma(x').$$
Since $\Gamma_0$ is flat, then the term $\int_{\Gamma_0}((\varphi\cdot\nabla){\bf n})\cdot \varphi\,d\sigma(x')$ vanishes, ans so we get (\ref{Gaffney}). Then, (\ref{Gaffney_div0}) is a consequence of the free divergence condition. The same result holds if $H_0^1(\Omega^\ep)^3$.

Finally, the change of variables (\ref{dilatacion}) applied to (\ref{Gaffney}) and (\ref{Gaffney_div0}) implies (\ref{Gaffney_tilde}) and (\ref{Gaffney_div0_tilde}), respectively, for any $\widetilde\varphi$ obtained from $\varphi$ by the change of variables (\ref{dilatacion}).
\end{proof}

\begin{lemma}For all $\widetilde\varphi \in \widetilde V^\ep$, the following inequalities hold\begin{equation}\label{trace_estimate}
 \|\widetilde\varphi\|_{L^2(\Gamma_0)^2}\leq \ep C\|D_\ep \widetilde\varphi\|_{L^2(\widetilde \Omega^\ep)^{3\times 3}}.
\end{equation}
\end{lemma}
\begin{proof}
Thank to $\widetilde \varphi(z',h(z'))=0$ in $\omega$, we have that
$$\begin{array}{rl}
\displaystyle \int_{\Gamma_0}|\widetilde \varphi|^2d\sigma= &\displaystyle \int_{\omega}|\widetilde \varphi(x',0)|^2\,dx'= \int_{\omega}\left|\int_0^{h(x')}\partial_{z_3}\widetilde \varphi(z',z_3)\,dz_3\right|^2\,dz'\leq   h_{\rm max}\int_{ \Omega}|\partial_{z_3}\widetilde \varphi|^2\,dz,
\end{array}$$
that is,
$$\|\widetilde \varphi'\|_{L^2(\Gamma_0)^2}\leq h_{\rm max}^{1\over 2}\|\partial_{z_3}\widetilde \varphi'\|_{L^2( \Omega)^{2}},$$
which implies (\ref{trace_estimate}).

\end{proof}

\subsection{Estimates for velocity and microrotation}

\begin{lemma}[Estimates for velocity and microrotation]\label{lemma_estimates}The following estimates hold for the solution $({\bf \widetilde u}_\ep, {\bf\widetilde w}_\ep)$ of problem (\ref{Form_Var_vel_tilde}):
\begin{equation} \label{estim_sol_dil1}
\|{\bf \widetilde u}_\epsilon\|_{L^2(\Omega)^3}\leq C\ep^2,  \quad
\|D_{\epsilon} {\bf  \widetilde u}_\epsilon\|_{L^2(\Omega)^{3\times 3}}\leq C\epsilon,\end{equation}
\begin{equation} \label{estim_sol_dil2}
\|{\bf \widetilde w}_\epsilon\|_{L^2(\Omega)^3}\leq C\ep, \quad
\|D_{\epsilon} {\bf \widetilde w}_\epsilon\|_{L^2(\Omega)^{3\times 3}}\leq C,
\end{equation}
\begin{equation}\label{uestimtrace}
\|{\bf \widetilde u}'_\ep\|_{L^2(\Gamma_0)^2}\leq C\ep^{3-\gamma\over 2}.
\end{equation}
\end{lemma}
\begin{proof}  From (\ref{Form_Var_vel_tilde}) with $\widetilde \varphi={\bf \widetilde u}_\ep$ and $\widetilde \psi={\bf \widetilde w}_\ep$, we have that
\begin{equation}\label{eq1}
\begin{array}{c}
\displaystyle
\|D_\ep{\bf \widetilde u}_\ep\|^2_{L^2(\Omega)^{3\times 3}}+\lambda\ep^{\gamma-1}\|{\bf\widetilde  u}'_\ep\|^2_{L^2(\Gamma_0)^2}=2N^2\int_{\Omega}{\rm rot}_\ep({\bf \widetilde w}_\ep)\cdot {\bf \widetilde u}_\ep\,dz+\int_{\Omega}{\bf f}'\cdot {\bf \widetilde u}_\ep\,dz,
\end{array}
\end{equation}
and
\begin{equation}\label{eq2}
\begin{array}{c}
\displaystyle
\ep^2R_c\|D_\ep {\bf \widetilde w}_\ep\|^2_{L^2(\Omega)^{3\times 3}}
+4N^2\|{\bf \widetilde w}_\ep\|^2_{L^2(\Omega)^3}=2N^2\int_{\Omega}{\rm rot}_\ep({\bf \widetilde u}_\ep)\cdot {\bf \widetilde w}_\ep\,dz+\ep\int_{\Omega}{\bf g}'\cdot{\bf \widetilde w}_\ep\,dz.
\end{array}
\end{equation}

Taking now into account that $\int_{\Omega}{\rm rot}_\ep({\bf \widetilde w}_\ep)\cdot {\bf \widetilde u}_\ep\,dz=\int_{\Omega}{\rm rot}_\ep({\bf \widetilde u}_\ep)\cdot {\bf \widetilde w}_\ep\,dz$, by Cauchy-Schwarz's inequality and (\ref{Gaffney_tilde}), we deduce that
\begin{equation}\label{eq3}
\begin{array}{c}
\displaystyle
2N^2\int_{\Omega}{\rm rot}_\ep({\bf \widetilde w}_\ep)\cdot {\bf \widetilde u}_\ep\,dz
\leq 2N^2\| {\rm rot}_\ep({\bf \widetilde u}_\ep)\|_{L^2(\Omega)^3}\|{\bf \widetilde w}_\ep\|_{L^2(\Omega)^3}\leq 2N^2\|D_\ep{\bf \widetilde u}_\ep\|_{L^2(\Omega)^{3\times 3}}\|{\bf \widetilde w}_\ep\|_{L^2(\Omega)^3},
\end{array}
\end{equation}
and by Cauchy-Schwarz's inequality and the estimate (\ref{Poincare2}), we have that
\begin{equation}\label{eq4}
\begin{array}{c}
\displaystyle
\int_{\Omega}{\bf \widetilde f}'\cdot {\bf \widetilde u}'_\ep\,dz
\leq \|{\bf f'}\|_{L^2(\Omega)^2}\|{\bf \widetilde u}_\ep\|_{L^2(\Omega)^3}
\leq h_{\rm max}^{1\over 2}C_p \ep \|{\bf f}'\|_{L^2(\omega)^2}\|D_\ep{\bf \widetilde u}_\ep\|_{L^2(\Omega)^{3\times 3}}.
\end{array}
\end{equation}
Taking into account (\ref{eq3})-(\ref{eq4}) in (\ref{eq1}), we deduce
\begin{equation}\label{eq1_1}
\begin{array}{c}
\displaystyle
\|D_\ep{\bf \widetilde u}_\ep\|_{L^2(\Omega)^{3\times 3}}\leq 2N^2\|{\bf \widetilde w}_\ep\|_{L^2(\Omega)^3}+ h_{\rm max}^{1\over 2}C_p \ep \|{\bf f}'\|_{L^2(\omega)^2}.
\end{array}
\end{equation}
On the other hand, by Cauchy-Schwarz's inequality and  (\ref{Gaffney_div0_tilde}), we deduce that
\begin{equation}\label{eq5}
\begin{array}{c}
\displaystyle
2N^2\int_{\Omega}{\rm rot}_\ep({\bf \widetilde u}_\ep)\cdot{\bf \widetilde w}_\ep\,dz
\leq 2N^2 \| D_\ep{\bf \widetilde u}_\ep\|_{L^2(\Omega)^{3\times 3}}\|{\bf \widetilde w}_\ep\|_{L^2(\Omega)^3},
\end{array}
\end{equation}
and by Cauchy-Schwarz's inequality and the estimate (\ref{Poincare2}), we have that
\begin{equation}\label{eq6}
\begin{array}{c}
\displaystyle
\ep\int_{\Omega}{\bf g}'\cdot{\bf \widetilde w}'_\ep\,dz
\leq h_{\rm max}^{1\over 2}\ep\|{\bf g'}\|_{L^2(\omega)^2}\|{\bf \widetilde w}_\ep\|_{L^2(\Omega)^3}.
\end{array}
\end{equation}
Taking into account (\ref{eq5})-(\ref{eq6}) in (\ref{eq2}), we deduce
\begin{equation}\label{eq2_1}
\begin{array}{c}
\displaystyle
2N^2\|{\bf \widetilde w}_\ep\|_{L^2(\Omega)^3}\leq N^2 \| D_\ep{\bf \widetilde u}_\ep\|_{L^2(\Omega)^{3\times 3}} + {1\over 2}h_{\rm max}^{1\over 2}\ep\|{\bf g'}\|_{L^2(\omega)^2}.
\end{array}
\end{equation}
From (\ref{eq1_1}) and (\ref{eq2_1}), we obtain
\begin{equation}\label{eq12_total}
\begin{array}{c}
\displaystyle
(1-N^2)\|D_\ep{\bf \widetilde u}_\ep\|_{L^2(\Omega)^{3\times 3}}\leq  {1\over 2}h_{\rm max}^{1\over 2}\ep\|{\bf g'}\|_{L^2(\omega)^2}+ h_{\rm max}^{1\over 2}C_p \ep \|{\bf f}'\|_{L^2(\omega)^2}.
\end{array}
\end{equation}
Taking into account the assumption of ${\bf f}'$ and ${\bf g}'$, the fact that $0<N^2<1$ and the inequality (\ref{Poincare2}), we obtain
\begin{equation}\label{eq12_total_integrating2}
\|D_\ep{\bf \widetilde u}_\ep\|_{L^2(\Omega)^{3\times 3}}\leq C\ep,\quad  \|{\bf \widetilde u}_\ep\|_{L^2(\Omega)^{3}}\leq C\ep^2.
\end{equation}
Now, let us prove the estimates for microrotation. To do this, by Cauchy-Schwarz's inequality, (\ref{Gaffney_div0_tilde}) and (\ref{Poincare2}), we deduce that
\begin{equation}\label{eq5222}
\begin{array}{c}
\displaystyle
2N^2\int_{\Omega}{\rm rot}_\ep({\bf \widetilde u}_\ep)\cdot{\bf \widetilde w}_\ep\,dz
\leq 2N^2 \| D_\ep{\bf \widetilde w}_\ep\|_{L^2(\Omega)^{3\times 3}}\|{\bf \widetilde u}_\ep\|_{L^2(\Omega)^3},
\end{array}
\end{equation}
and
\begin{equation}\label{eq52221}
\begin{array}{c}
\displaystyle
\ep\int_{\Omega}{\bf g}'\cdot{\bf \widetilde w}'_\ep\,dz
\leq h_{\rm max}^{1\over 2}\ep^2\|{\bf g'}\|_{L^2(\omega)^2}\|D_\ep{\bf \widetilde w}_\ep\|_{L^2(\Omega)^{3\times 3}}.
\end{array}
\end{equation}
 Taking into account (\ref{eq5222})-(\ref{eq52221}) in (\ref{eq2}), we deduce
 \begin{equation}\label{eq52223}
  \ep^2R_c\|D_\ep {\bf \widetilde w}_\ep\|_{L^2(\Omega)^{3\times 3}}
\leq  2N^2 \|{\bf \widetilde u}_\ep\|_{L^2(\Omega)^3}+h_{\rm max}^{1\over 2}\ep^2\|{\bf g'}\|_{L^2(\omega)^2},
\end{equation}
and from estimates of ${\bf \widetilde u}_\ep$ and assumption on ${\bf g}'$, we deduce
 \begin{equation}\label{eq52224}
  \ep^2R_c\|D_\ep {\bf \widetilde w}_\ep\|_{L^2(\Omega)^{3\times 3}}
\leq  C\ep^2,
\end{equation}
which implies (\ref{estim_sol_dil2}).

Finally, let us derive the estimates for ${\bf \widetilde u}_\ep'$ on $\Gamma_0$. To do this, from (\ref{eq1}), (\ref{eq3}) and (\ref{eq4}) and estimates for velocity and microrotation, we deduce
\begin{equation}\label{eq1trace}
\begin{array}{c}
\displaystyle
\lambda\ep^{\gamma-1}\|{\bf\widetilde  u}'_\ep\|^2_{L^2(\Gamma_0)^2}\leq C\ep^2,
\end{array}
\end{equation}
which implies (\ref{uestimtrace}).
\end{proof}

\subsection{Estimates for pressure}
The estimates for the pressure will be obtained by using a decomposition result for the pressure $p_\ep$ in $\Omega^\ep$.
\begin{proposition}
The pressure $p_\ep\in L^2_0(\Omega^\ep)$ satisfies the following decomposition
\begin{equation}\label{decomp}
p_\ep=p_\ep^{(0)}+p_\ep^{(1)},
\end{equation}
where $p_\ep^{(0)}\in H^1(\omega)$, independent of $x_3$, and $p_\ep^{(1)}\in L^2(\Omega^\ep)$.  Moreover, the following estimates hold
\begin{equation}\label{estimdecomp}
\|p_\ep^{(0)}\|_{H^1(\omega)}\leq C\ep^{-{3\over 2}}\|\nabla p_\ep\|_{H^{-1}(\Omega^\ep)^3},\quad \|p_\ep^{(1)}\|_{H^1(\omega)}\leq C\|\nabla p_\ep\|_{H^{-1}(\Omega^\ep)^3}.
\end{equation}
\end{proposition}
\begin{proof}
This result follows from  \cite[Corollary 3.4]{CLS} by assuming
$$\Omega^\ep=\omega\times (\varepsilon \vartheta),\quad \forall\varepsilon>0,$$ with $\omega\subset \mathbb{R}^2, \vartheta\subset \mathbb{R}$ connected Lipschitz open sets and exponent $q=2$, see \cite[Remark 3.1]{CLS}. In the present case, $\vartheta=(0,h(x'))$ with $h\in W^{1,\infty}$ as we set in Section \ref{sec:setting}. Thus, $\Omega^\ep$ satisfies the abstract conditions (3.3)-(3.5) given in \cite{CLS}, and so,   \cite[Corollary 3.4]{CLS} can be applied.
\end{proof}
\begin{corollary}\label{cor:pressure}
Denoting by $\widetilde p_\ep^{(1)}$ the rescaled function associated with $p_\ep^{(1)}$, then  the following decomposition holds
\begin{equation}\label{decomptildeestim}
\widetilde p_\ep=p_\ep^{(0)}+\widetilde p_\ep^{(1)},
\end{equation}
where the pressures $p_\ep^{(0)}$  and $\widetilde p_\ep^{(1)}$ satisfy the following estimates
\begin{equation}\label{estimpressure}
\|p_\ep^{(0)}\|_{H^1(\omega)}\leq C,\quad \|\widetilde p_\ep^{(1)}\|_{L^2(\Omega^\ep)}\leq C\ep.
\end{equation}
\end{corollary}

\begin{proof}  From estimates (\ref{estimdecomp}) and taking into account that
$$\|\widetilde p_\ep^{(1)}\|_{L^2(\Omega)}=\ep^{-{1\over 2}}\| p_\ep^{(1)}\|_{L^2(\Omega^\ep)},$$
then to derive (\ref{estimpressure}), we just need to derive the following estimate
\begin{equation}\label{nablap}
\|\nabla p_\ep\|_{H^{-1}(\Omega^\ep)^3}=\sup_{\varphi\in H^1_0(\Omega^\ep)^3}
{|\langle \nabla p_\ep,\varphi\rangle_{H^{-1}(\Omega^\ep)^3,H_0^1(\Omega^\ep)^3}|\over \|\varphi\|_{H^1_0(\Omega^\ep)^3}}\leq C\ep^{3\over 2}.
\end{equation}
To to this, we take into account that, according to the change of variables (\ref{dilatacion}) and estimates for velocity ${\bf \widetilde u}_\ep$ and ${\bf \widetilde w}_\ep$ given in Lemma \ref{lemma_estimates}, we have
$$\|{\bf u}_\ep\|_{L^2(\Omega^\ep)^3}\leq C\ep^{5\over 2},\quad \|D{\bf u}_\ep\|_{L^2(\Omega^\ep)^{3\times 3}}\leq C\ep^{3\over 2}, \quad\|{\bf w}_\ep\|_{L^2(\Omega^\ep)^3}\leq C\ep^{3\over 2},\quad \|D{\bf w}_\ep\|_{L^2(\Omega^\ep)^{3\times 3}}\leq C\ep^{1\over 2}.
$$
Now, taking $\varphi\in H^1_0(\Omega^\ep)^3$ as test function in (\ref{Form_Var_vel}), using the Cauchy-Schwarz inequality,  the inequalities (\ref{Poincare}) and (\ref{Gaffney}) and the previous estimates of velocity and microrotation, we get
$$\begin{array}{rl}\displaystyle |\langle \nabla p_\ep,\varphi\rangle_{H^{-1}(\Omega^\ep)^3,H_0^1(\Omega^\ep)^3}|\leq &\displaystyle\|D {\bf u}_\ep\|_{L^2(\Omega^\ep)^{3\times 3}}\|D\varphi\|_{L^2(\Omega^\ep)^{3\times 3}}+2N^2\|{\rm rot} {\bf w}_\ep\|_{L^2(\Omega^\ep)^{3}}\|\varphi\|_{L^2(\Omega^\ep)^{3}}\\
\noame
&\displaystyle+h_{\rm max}^{1\over 2}\ep^{1\over 2} \|{\bf f'}\|_{L^2(\omega)^3}\|\varphi\|_{L^2(\Omega^\ep)^{3}}\\
\noame
\leq &\displaystyle\|D {\bf u}_\ep\|_{L^2(\Omega^\ep)^{3\times 3}}\|D\varphi\|_{L^2(\Omega^\ep)^{3\times 3}}+2N^2\ep\|D {\bf w}_\ep\|_{L^2(\Omega^\ep)^{3\times 3}}\|D\varphi\|_{L^2(\Omega^\ep)^{3\times 3}}\\
\noame
&\displaystyle+h_{\rm max}^{1\over 2}\ep^{3\over 2} \|{\bf f'}\|_{L^2(\omega)^3}\|D\varphi\|_{L^2(\Omega^\ep)^{3\times 3}}\\
\noame
 \leq &\displaystyle C\ep^{3\over 2}\|\varphi\|_{H^1_0(\Omega^\ep)^3},
\end{array}$$
which implies (\ref{nablap}).
\end{proof}

\subsection{Convergences}
In this section, depending on the value of $\gamma$, we give  convergence results concerning the asymptotic behavior of the sequence $({\bf \widetilde u}_\ep, {\bf \widetilde w}_\ep, \widetilde p_\ep)$  satisfying the {\it a priori} estimates given in Lemmas \ref{lemma_estimates} and Corollary \ref{cor:pressure}.  To do this, let us introduce the following Hilbert space
$$V_{z_3}(\Omega)=\{\widetilde \varphi\in L^2(\Omega)\, :\, \partial_{z_3}\varphi\in L^2(\Omega)\},$$
with the norm
$$\|\widetilde\varphi\|_{V_{z_3}}^2=\|\widetilde\varphi\|^2_{L^2(\Omega)}+\|\partial_{z_3}\widetilde\varphi\|^2_{L^2(\Omega)}.$$
\begin{lemma}\label{lem_asymp_sub}
For a subsequence of $\ep$ still denote by $\ep$, we have the following convergence results:
\begin{itemize}
\item[(i)] {\it (Velocity)} There exists ${\bf \widetilde u}=({\bf \widetilde u'},\widetilde u_{3})\in V_{z_3}(\Omega)^3$, with $\widetilde u_{3}\equiv 0$, such that
\begin{eqnarray}
&\displaystyle \ep^{-2}{\bf  \widetilde u}^\ep\rightharpoonup {\bf \widetilde u}\quad \hbox{in}\quad V_{z_3}(\Omega)^3,\label{conv_u_sub_tilde}\\
\noame
&\displaystyle {\rm div}_{x'}\left(\int_0^{h(z')}{\bf \widetilde u}'(z',z_3)\,dz_3\right)=0\quad \hbox{  in  }\omega, &\label{div_x_sub_tilde}\\
\noame
&\displaystyle \left(\int_0^{h(z')}{\bf \widetilde u}'(z',z_3)\,dz_3\right)\cdot {\bf n}'=0\quad \hbox{  on  }\partial\omega,&\nonumber
 \end{eqnarray}
 where ${\bf n}'$ is the outward unit normal vector to the boundary $\omega$. Moreover, the following values  of the velocity on the boundaries hold:
 \begin{itemize}
 \item If $\gamma\geq -1$, then ${\bf \widetilde u}'$ satisfies the no-slip condition on the top boundary, i.e.
 \begin{equation}\label{case1bc}
 {\bf \widetilde u}'=0\quad\hbox{on}\quad \Gamma_1.
 \end{equation}
  \item If $\gamma< -1$, then ${\bf \widetilde u}'$ satisfies the no-slip condition on the top and bottom boundaries, i.e.
 \begin{equation}\label{case2bc}
 {\bf \widetilde u}'=0\quad\hbox{on}\quad \Gamma_0\cup \Gamma_1.
 \end{equation}
 \end{itemize}
\item[(ii)] {\it (Microrotation)} There exist ${\bf \widetilde w}=({\bf \widetilde w}',\widetilde w_3)\in V_{z_3}(\Omega)^3$, with $\widetilde w_3\equiv 0$ and ${\bf \widetilde w}'=0$ on $ \Gamma_0\cup \Gamma_1$, such that
\begin{eqnarray}
&\displaystyle \ep^{-1}{\bf \widetilde w}_\ep\rightharpoonup {\bf\widetilde w}\quad \hbox{in}\quad V_{z_3}(\Omega)^3.\label{conv_w_sub_tilde}
\end{eqnarray}
\item[(iii)] {\it (Pressure)} There  exist a function  $\widetilde p\in L^2_0(\Omega)\cap H^1(\omega)$, i.e. $\widetilde p$ is independent of $z_3$ and has mean value zero in $\Omega$, and a function $\widetilde p^{(1)}\in L^2(\Omega)$  such that
\begin{eqnarray}
&\displaystyle  \widetilde p_\ep\to  \widetilde p\quad\hbox{in}\quad L^2(\Omega),\quad  \ep^{-1}\partial_{z_3}\widetilde p_\ep\to  \partial_{z_3}\widetilde p^{(1)}\quad\hbox{in}\quad L^2(\omega;H^{-1}(0,h(z'))).& \label{conv_P_sub}
\end{eqnarray}
\end{itemize}
\end{lemma}
\begin{proof} We start proving $(i)$. We will only give some remarks and,  for more details,  we refer the reader to in \cite{BayadaLuc, Bayada_Gamouana, Bayada_Gamouana2, grau1}.  The proof of (\ref{conv_u_sub_tilde}), with $\widetilde u_3\equiv 0$ and ${\bf \widetilde u}'=0$ on $\Gamma_1$, and (\ref{div_x_sub_tilde}) is standard, so we focus on the derivation of the boundary conditions satisfied by ${\bf \widetilde u}'$ depending on the value of $\gamma$. From the trace estimate (\ref{trace_estimate}), we deduce from (\ref{estim_sol_dil1}) that
$$\|{\bf \widetilde u}'_\ep\|_{L^2(\Gamma_0)^2}\leq C\ep\|D_\ep {\bf \widetilde u}_\ep\|_{L^2(\Omega)^{3\times 3}}\leq C\ep^2,$$
which implies
\begin{equation}\label{convgamma0}\ep^{-2}{\bf \widetilde u}_\ep\rightharpoonup {\bf\widetilde u}'\quad\hbox{in}\quad L^2(\Gamma_0)^2.
\end{equation}
Moreover,  from estimate (\ref{uestimtrace}), we deduce
\begin{equation}\label{estimgamma0}
\|{\bf \widetilde u}'_\ep\|_{L^2(\Gamma_0)^2}\leq C\min\{\ep^{3-\gamma\over 2}, \ep^2\}.
\end{equation}
Thus, there is a critical value when ${3-\gamma\over 2}=2$, i.e. when $\gamma=-1$. Depending on this value, we deduce the following:
\begin{itemize}
\item In the case $\gamma\geq -1$, then the optimal quantity in (\ref{estimgamma0}) is $C\ep^2$. Then, we still have the convergence (\ref{convgamma0}).
\item In the case $\gamma<-1$, then the optimal quantity in (\ref{estimgamma0}) is $C\ep^{3-\gamma\over 2}$, and since $\ep^{-{\gamma+1\over 2}}\to 0$ and $\ep^{\gamma-3\over 2}{\bf \widetilde u'}_\ep$ is bounded in $L^2(\Gamma_0)^2$, we have
$$\ep^{-2}{\bf \widetilde u}_\ep=\ep^{-{\gamma+1\over 2}}(\ep^{\gamma-3\over 2}{\bf\widetilde u}'_\ep)\rightharpoonup 0\quad\hbox{in}\quad L^2(\Gamma_0)^2.$$
Then, from the uniqueness of the limit, ${\bf \widetilde u}'=0$ on $\Gamma_0$.

\end{itemize}

We continue proving $(ii)$.  According to estimates (\ref{estim_sol_dil2}), the proof of convergence (\ref{conv_w_sub_tilde}) and the boundary condition ${\bf \widetilde w}=0$ on $\Gamma_0\cup \Gamma_1$ is standard, see for instance \cite{BayadaLuc, Bayada_Gamouana, Bayada_Gamouana2,  SG1}.   It remains to prove that $\widetilde w_3\equiv0$. To do this, we consider as test function $\widetilde \psi(z)=\ep^{-1}(0,0,\widetilde \psi_3)$ with $\widetilde\psi_3\in \mathcal{D}(\Omega)$,  in the variational formulation (\ref{Form_Var_vel_tilde}), and so we get
$$
\begin{array}{l}
\displaystyle \ep R_c\int_{\Omega}\nabla_{x'}\widetilde w_{\ep,3}\cdot\nabla_{x'}\widetilde \psi_3\,dz+\ep^{-1} R_c\int_{\Omega}\partial_{z_3}\widetilde w_{\ep,3}\,\partial_{z_3}\widetilde \psi_3\,dz+4N^2\ep^{-1}\int_{\Omega}{\widetilde w}_{\ep,3}\, \widetilde \psi_3\,dz-2N^2\int_\Omega \ep^{-1}{\rm Rot}_{z'}({\widetilde u}'_\ep)\,\widetilde\psi_3\,dz=0.
\end{array}
$$
Passing to the limit by using convergences of  ${\bf \widetilde u}_\ep$ and ${\bf \widetilde w}_\ep$ given respectively in (\ref{conv_w_sub_tilde}), we get
$$R_c\int_\Omega \partial_{z_3}\widetilde w_3\,\partial_{z_3}\widetilde \psi_3\,dz+4N^2\int_\Omega \widetilde w_3\,\widetilde \psi_3\,dz=0\,,$$
and taking into account that $\widetilde w_3=0$ on $\Gamma_1\cup \Gamma_0$,  it is easily deduced that $\widetilde w_3\equiv 0$  in $\Omega$.\\

We finish the proof with $(iii)$. Estimate (\ref{estimpressure}) implies, up to a subsequence, the existence of $\widetilde p\in H^1(\omega)$ and $\widetilde p^{(1)}\in L^2(\Omega)$ such that
\begin{equation}\label{conv_P_sub}
 \widetilde p^{\widetilde (0)}_\ep\rightharpoonup \widetilde p\quad \hbox{  in  }H^1(\omega),\quad \ep^{-1}\widetilde p^{(1)}_\ep\rightharpoonup \widetilde p^{(1)}\quad \hbox{  in  }L^2(\Omega).
\end{equation}
As consequece,  by  the decomposition of $\widetilde p_\ep$ given in (\ref{decomptildeestim}), we deduce convergence (\ref{conv_P_sub}).

Finally, since $\widetilde p_\ep$ has mean value zero in $\Omega$ and $\widetilde p$ does not depend on $z_3$, it holds
$$
  0=\int_\Omega  \widetilde p_\ep\,dz\to \int_\Omega\widetilde p\,dz= \int_{\omega}h(z')\widetilde p(z')\,dz'=0,
$$
which implies $\widetilde p\in L^2_0(\Omega)$.
\end{proof}

\section{Limit problem}\label{sec:limitproblem}

\noindent Using convergences given in Lemma \ref{lem_asymp_sub}, in this section we give the reduced  effective system satisfied by $({\bf \widetilde u},{\bf \widetilde w}, \widetilde p)$ with boundary conditions on $\Gamma_0$ depending on the value of $\gamma$.

\begin{theorem}[Limit problem]\label{thm_sub1}
The triplet  of functions
   $({\bf \widetilde u}, {\bf \widetilde w}, \widetilde p)\in V_{z_3}(\Omega)^3\times V_{z_3}(\Omega)^3 \times  (L^2_0(\Omega)\cap H^1(\omega))$, with $\widetilde u_3=\widetilde w_3\equiv 0$, given in Lemma \ref{lem_asymp_sub} is the unique solution of the reduced effective micropolar system
   \begin{equation}\label{hom_system_sub_u}
 \begin{array}{rcll}
\displaystyle
-\partial_{z_3}^2 {\bf \widetilde u}'(z) -2N^2{\rm rot}_{z_3}({\bf \widetilde  w}'(z))&=&\displaystyle  {\bf f}'(z')- \nabla_{z'}\widetilde p(z')&\hbox{ in }\Omega,\\
\noame
\displaystyle
-R_c\partial_{z_3}^2   {\bf \widetilde w}'(z)+4N^2 {\bf \widetilde  w}'(z) -2N^2{\rm rot}_{z_3}({\bf \widetilde u'}(z))&=&{\bf g}'(z')\displaystyle &\hbox{ in }\Omega,\\
\end{array} 
\end{equation}
with divergence condition
\begin{eqnarray}
\displaystyle -{\rm div}_{z'}\left( \int_0^{h(z')}{\bf \widetilde u}'(z)\,dz_3\right)=0&&\hbox{ in }\omega,\label{divzh_limit}
\end{eqnarray}
the no-slip condition for velocity and no-spin condition for microrotation on the top boundary
\begin{eqnarray}
{\bf \widetilde u}'=0,\quad{\bf \widetilde w}'=0&&\hbox{on } \Gamma_1,\label{bc_limit1}
\end{eqnarray}
and the following boundary conditions on the bottom boundary $\Gamma_0$ depending on the value of $\gamma$:
\begin{itemize}
\item If $\gamma>-1$, a perfect slip boundary for velocity and no-spin condition for microrotation
\begin{equation}\label{gmayor1}
-\partial_{z_3}{\bf \widetilde u}'=0, \quad {\bf \widetilde w}'=0\quad\hbox{on}\quad \Gamma_0.
\end{equation}
\item If $\gamma=-1$, a partial slip boundary for velocity and no-spin condition for microrotation
\begin{equation}\label{gigual1}
-\partial_{z_3}{\bf \widetilde u}'=-\lambda {\bf u}', \quad {\bf \widetilde w}'=0\quad\hbox{on}\quad \Gamma_0.
\end{equation}
\item If $\gamma<-1$, a no-slip boundary for velocity and no-spin condition for microrotation
\begin{equation}\label{gmenor1}
{\bf \widetilde u}'=0, \quad {\bf \widetilde w}'=0\quad\hbox{on}\quad \Gamma_0.
\end{equation}
\end{itemize}
\end{theorem}
\begin{proof} From Lemma \ref{lem_asymp_sub}, it remains to prove (\ref{hom_system_sub_u}) and boundary conditions for velocity on the bottom $\Gamma_0$. We divide the proof in three steps.  \\

\noindent {\it Step 1.} We prove  (\ref{hom_system_sub_u})$_{1,2}$ with boundary condition on $\Gamma_0$ depending on the value of $\gamma$. According to Lemma \ref{lem_asymp_sub}, we consider  in (\ref{Form_Var_vel_tilde})$_1$ where $\widetilde\varphi \in C^1(0,h(z');C_c^1(\omega)^3)$ with $\widetilde \varphi_3\equiv 0$, $\widetilde \varphi'=0$ on $\Gamma_1$ .  Moreover, if $\gamma<-1$, we also consider $\widetilde\varphi'=0$ on $\Gamma_0$. This gives
\begin{equation}\label{Form_Var_vel_tilde_hat_v2_proof1}
\begin{array}{l}
\displaystyle  \int_{\Omega}D_{z'}{\bf\widetilde  u}'_\ep: D_{z'}\widetilde\varphi'\,dz+\ep^{-2}\int_{\Omega}\partial_{z_3}{\bf\widetilde  u}'_\ep \cdot \partial_{z_3}\widetilde\varphi'\,dz+\lambda\ep^{\gamma-1}\int_{\Gamma_0}{\bf \widetilde u}_\ep'\cdot \widetilde\varphi'\,d\sigma(z')\\
\noame
\displaystyle - \int_{\Omega}\widetilde p_\ep\,{\rm div}_{z'}(\widetilde\varphi')\,dz-2N^2\ep^{-1}\int_{\Omega} {\rm rot}_{z_3}({\bf \widetilde w}'_\ep)\cdot \widetilde\varphi'\,dz =\int_{\Omega}{\bf f}'\cdot \widetilde\varphi'\,dz.
\end{array}
\end{equation}
Below, let us pass to the limit when $\epsilon$ tends to zero in each term of the previous variational formulation:
\begin{itemize}
\item First and second terms. Using estimate of $D_{z'}{\bf \widetilde u}_\ep$ given in (\ref{estim_sol_dil1}) and convergence (\ref{conv_u_sub_tilde}), we get
$$\begin{array}{rl}
\displaystyle \int_{\Omega}D_{z'}{\bf\widetilde  u}'_\ep: D_{z'}\widetilde\varphi'\,dz+\ep^{-2}\int_{\Omega}\partial_{z_3}{\bf\widetilde  u}'_\ep \cdot \partial_{z_3}\widetilde\varphi'\,dz=&\displaystyle
\int_{\Omega}\partial_{z_3}{\bf\widetilde  u}'\cdot \partial_{z_3}\widetilde\varphi'\,dz+O_\ep.
\end{array}
$$
\item Third term.  Using convergence of ${\bf \widetilde u}'_\ep$ in $\Gamma_0$ given in (\ref{convgamma0}) depending on $\gamma$:
\begin{itemize}
\item If $\gamma=-1$, then
$$\lambda\ep^{-2}\int_{\Gamma_0}{\bf \widetilde u}_\ep'\cdot \widetilde\varphi\,d\sigma(z')=\lambda \int_{\Gamma_0}{\bf \widetilde u}'\cdot \widetilde\varphi\,d\sigma(z')+O_\ep.$$
\item If $\gamma>-1$, since $\ep^{\gamma+1}\to 0$, then we deduce
$$
\displaystyle \lambda\ep^{\gamma-1}\int_{\Gamma_0}{\bf \widetilde u}_\ep'\cdot \widetilde\varphi\,d\sigma(z')=\displaystyle
\lambda\ep^{\gamma+1}\int_{\Gamma_0}\ep^{-2}{\bf \widetilde u}_\ep'\cdot \widetilde\varphi\,d\sigma(z')\to 0.$$
\item If $\gamma<-1$, since $\varphi'=0$ on $\Gamma_0$ in this case, we have
\begin{equation}\label{case3_zero}\lambda\ep^{\gamma-1}\int_{\Gamma_0}{\bf \widetilde u}_\ep'\cdot \widetilde\varphi\,d\sigma(z')=0.
\end{equation}
\end{itemize}
\item Fourth term. Using convergence (\ref{conv_P_sub}) of the pressure $\widetilde p_\ep$, we get
$$-\int_{\Omega}\widetilde p_\ep\,{\rm div}_{z'}(\widetilde\varphi')\,dz= -\int_{\Omega}\widetilde p\,{\rm div}_{z'}(\widetilde\varphi')\,dz+O_\ep.$$
\item Fifth term. Using convergence (\ref{conv_w_sub_tilde}) and integration by parts, we get
$$
-2N^2\ep^{-1}\int_{\Omega} {\rm rot}_{z_3}({\bf \widetilde w}'_\ep)\cdot \widetilde\varphi'\,dz
=-2N^2\int_{\Omega} {\rm rot}_{z_3}({\bf \widetilde w}')\cdot \widetilde\varphi'\,dz+O_\ep.
$$
\end{itemize}
Therefore, by previous convergences, we deduce the following general limit variational formulation
\begin{equation}\label{form_var_limit_vel}\int_{\Omega}\partial_{z_3}{\bf\widetilde  u}'\cdot \partial_{z_3}\widetilde\varphi'\,dz+\lambda \int_{\Gamma_0}{\bf \widetilde u}'\cdot \widetilde\varphi\,d\sigma(z')-\int_{\Omega}\widetilde p\,{\rm div}_{z'}(\widetilde\varphi')\,dz-2N^2\int_{\Omega} {\rm rot}_{z_3}({\bf \widetilde w}')\cdot \widetilde\varphi'\,dz=\int_\Omega {\bf f}'\cdot \widetilde\varphi'\,dz, \end{equation}
where, by density, 
\begin{itemize}
\item If $\gamma\geq -1$, then $\widetilde\varphi'\in V_{z_3}(\Omega)^2$ with $\widetilde\varphi'=0$ on $\Gamma_1$. Moreover, if $\gamma=-1$ then we can take in (\ref{form_var_limit_vel}) the value $\lambda\in (0,+\infty)$, and if $\gamma>-1$ then we can take in (\ref{form_var_limit_vel}) the value $\lambda=0$ because the boundary term disappear.
\item If $\gamma<-1$, then $\varphi'\in V_{z_3}(\Omega)^2$ with $\varphi'=0$ on $\Gamma_0\cup \Gamma_1$ and we can take $\lambda=+\infty$ so the boundary term disappear (see (\ref{case3_zero})).
\end{itemize}
Integrating by parts  and  the definition of the rotational given in (\ref{rotationals}), then (\ref{form_var_limit_vel}) is equivalent to problem (\ref{hom_system_sub_u})$_{1,2}$ with boundary conditions on the bottom (\ref{gmayor1})--(\ref{gmenor1}) depending on the value of $\gamma$.\\

\noindent {\it Step 2.} We prove   (\ref{hom_system_sub_u})$_{3,4}$. According to Lemma \ref{lem_asymp_sub}, we consider  in (\ref{Form_Var_vel_tilde})$_2$ with $\widetilde \psi \in C^1_c(\Omega)^3$ with $\widetilde \psi_3\equiv 0$ in $\Omega$. After multiplying by $\ep^{-1}$, this gives
\begin{equation}\label{Form_Var_vel_tilde_hat_w2_proof1}
\begin{array}{l}
\displaystyle
\ep R_c\int_{\Omega}D_{z'}{\bf \widetilde w}'_\ep: D_{z'}\widetilde\psi'\,dz+ R_c\int_{\Omega}\ep^{-1}\partial_{z_3}{\bf \widetilde w}'_\ep\cdot \partial_{z_3}\widetilde\psi'\,dz+4N^2\int_{\Omega}\ep^{-1}{\bf \widetilde w}_\ep'\cdot\widetilde\psi'\,dz\\
\noame
\displaystyle -2N^2\int_{\Omega}\ep^{-2}{\rm rot}_{z_3}({\bf\widetilde u}'_\ep)\cdot \widetilde\psi'\,dz=\int_{\Omega}{\bf g}'\cdot \widetilde\psi'\,dz.
\end{array}
\end{equation}
 Below, let us pass to the limit when $\epsilon$ tends to zero in each term of the previous variational formulation:
\begin{itemize}
\item First and second term. Using estimate for $D_{z'}{\bf \widetilde w}'_\ep$ given in (\ref{estim_sol_dil2}) and convergence (\ref{conv_w_sub_tilde}), we get
$$\begin{array}{rl}
\displaystyle \ep R_c\int_{\Omega}D_{z'}{\bf \widetilde w}'_\ep: D_{z'}\widetilde\psi'\,dz+ R_c\int_{\Omega}\ep^{-1}\partial_{z_3}{\bf \widetilde w}'_\ep\cdot \partial_{z_3}\widetilde\psi'\,dz=&\displaystyle
 R_c\int_{\Omega}\partial_{z_3}{\bf \widetilde w}'\cdot \partial_{z_3}\widetilde \psi'\,dz+O_\ep.
\end{array}$$
\item Third term. Using convergence (\ref{conv_w_sub_tilde}), we get
$$\begin{array}{rl}
\displaystyle 4N^2\int_{\Omega}\ep^{-1}{\bf \widetilde w}_\ep'\cdot\widetilde\psi'\,dz=&\displaystyle
4N^2\int_{\Omega} {\bf \widetilde w}'\cdot\widetilde\psi'\,dz+O_\ep.
\end{array}$$
\item Fourth term. Using convergence (\ref{conv_u_sub_tilde}), we get
$$\begin{array}{rl}
\displaystyle
-2N^2\int_{\Omega}\ep^{-2}{\rm rot}_{z_3}({\bf\widetilde u}'_\ep)\cdot \widetilde\psi'\,dz=&\displaystyle
-2N^2\int_{\Omega}{\rm rot}_{z_3}({\bf\widetilde u}')\cdot \widetilde\psi'\,dz+O_\ep.
\end{array}$$
\end{itemize}
Therefore, by previous convergences, we deduce that the limit variational formulation is given by
\begin{equation}\label{form_var_1_changvar_hat1_limit2_sub2_w}
 \begin{array}{l}
\displaystyle  R_c\int_{\Omega}\partial_{z_3}{\bf \widetilde w}'\cdot \partial_{z_3}\widetilde \psi'\,dz4N^2\int_{\Omega} {\bf \widetilde w}'\cdot\widetilde\psi'\,dz-2N^2\int_{\Omega}{\rm rot}_{z_3}({\bf\widetilde u}')\cdot \widetilde\psi'\,dz=\int_{\Omega}{\bf g}'\cdot \widetilde\psi'\,dz,
\end{array}\end{equation}
for every $\widetilde \psi'\in V_{z_3}(\Omega)^2$. Integrating by parts and taking into account the definition of the rotational given in (\ref{rotationals}),  this variational formulation is equivalent to problem (\ref{hom_system_sub_u})$_{3,4}$ with no-spin condition on $\Gamma_0\cup\Gamma_1$.\\

{\it Step 3. Conclusion. }   To ensure that the whole sequence $(\ep^{-2}{\bf \widetilde u}_\ep, \ep^{-1} {\bf \widetilde w}_\ep,  \widetilde p_\ep)$ converges, it remains to prove the existence and uniqueness of weak solution of the effective system  (\ref{hom_system_sub_u}). This follows the lines of the proof of Theorem \ref{thm:existence}, so we omit it.

\end{proof}

\section{Expressions of velocity and microrotation}\label{sec:expressions}
 First, we give the expressions of ${\bf \widetilde u}$ and ${\bf \widetilde w}$ by solving problem (\ref{hom_system_sub_u}) with boundary condition (\ref{bc_limit1}) on $\Gamma_1$ and the different boundary conditions   (partial slip, no-slip and perfect slip) on $\Gamma_0$.

\begin{lemma}[Partial slip case]\label{expression_uw_hat_alpha}
The solution $({\bf \widetilde u}', {\bf \widetilde w}')$ of (\ref{hom_system_sub_u}) with boundary conditions (\ref{bc_limit1}) on $\Gamma_1$ and (\ref{gigual1}) on $\Gamma_0$ with $\lambda>0$ are given by the following expressions:
\begin{equation}\label{partial_u_expression}
  \begin{array}{l}
 {\bf  \widetilde u}'(z) ={1\over 2(1-N^2)}(\nabla_{z'}\widetilde p(z')-{\bf f}'(z'))\left[z_3^2-{h(z')}\left({2N^2\over k}\sinh(kz_3)-2z_3\right)A^\lambda_1(z')-{2N^2\over k} {h(z')}\left(\cosh(kz_3)-1\right)B^\lambda_1(z')
  \right]\\
  \noame
 +{1\over 2N^2} ({\bf g}'(z'))^\perp\left[\left(z_3-{1\over \lambda}\right)-\left(h(z') \right) \left({2N^2\over k}\sinh(kz_3)-2z_3\right)A^\lambda_2(z')-{2N^2\over k}\left(h(z')\right)\left(\cosh(kz_3)-1\right) B^\lambda_2(z')\right] \\
 \noame
 +{1\over \lambda}\left[-{h(z')}(\nabla_{z'}\widetilde p(z')-{\bf f}'(z'))A^\lambda_1(z') \right.
 \\
 \noame
 \qquad\left. +{1\over 2N^2}({\bf g}'(z'))^\perp
 \left(-1+\left(-2h(z')(1-N^2)+{2N^2\over k}\sinh(kz_3)-2z_3\right)A^\lambda_2(z')+{2N^2\over k}(\cosh(kz_3)-1)B^\lambda_2(z')
 \right)\right]\\
 \noame
 +{1\over \lambda^2}{1-N^2\over N^2}({\bf g}'(z')^\perp A^\lambda_2(z'),
 \end{array}
 \end{equation}
$$\begin{array}{l}
 {\bf \widetilde w}'(z)={1\over 2(1-N^2)}(\nabla_{z'}\widetilde p(z')-{\bf f}'(z'))^\perp\Big[z_3-h(z') (\cosh(kz_3)-1)A^\lambda_1(z')-h(z')\sinh(kz_3)B_1^\lambda(z')
  \Big]\\
  \noame
  +{h(z')\over 2N^2}{\bf g}'(z')\Big[  (\cosh(kz_3)-1)A^\lambda_2(z') +  \sinh(kz_3)B^\lambda_2(z')
  \Big]\\
  \noame
  -{1\over \lambda}{1\over 2N^2}{\bf g}'(z')\Big[(\cosh(kz_3)-1)A^\lambda_2(z')+\sinh(kz_3)B^\lambda_2(z')\Big],
 \end{array}
$$
where
\begin{equation}\label{parameterk}
k=2N\sqrt{{1-N^2\over R_c}},
\end{equation}
and  $A^\lambda_i, B^\lambda_2$, $i=1,2$, are given by
\begin{equation}\label{AB_partial_expression}\begin{array}{rl}
\noame
A^\lambda_1(z')=&\displaystyle -{1\over 2}{h(z')-{2N^2\over k}\tanh\left({kh(z')\over 2}\right)\over h(z')-{2N^2\over k}\tanh\left({kh(z')\over 2}\right)-{1\over \lambda}(1-N^2)},\\
\noame
B^\lambda_1(z')=&\displaystyle {1\over 2}\coth\left({kh(z')\over 2}\right){{2N^2\over k}\tanh\left({kh(z')\over 2}\right)-h(z')+{2\over\lambda}(1-N^2){1 \over 1+\cosh(kh(z'))}
\over {2N^2\over k}\tanh\left({kh(z')\over 2}\right)-h(z')+{1\over \lambda}(1-N^2)},\\
\noame
A^\lambda_2(z')=&\displaystyle{1\over 2}{1\over  {2N^2\over k}\tanh\left({kh(z')\over 2}\right)-h(z')+{1\over \lambda}(1-N^2)},\\
\noame
B^\lambda_2(z')&\displaystyle-{1\over 2}{\tanh\left({kh(z')\over 2}\right) \over  {2N^2\over k}\tanh\left({kh(z')\over 2}\right)-h(z')+{1\over \lambda}(1-N^2)}.
\end{array}
\end{equation}
\end{lemma}
\begin{proof} We have to solve the system (\ref{hom_system_sub_u}), i.e.
$$
\left\{\begin{array}{rcll}
\displaystyle
-\partial_{z_3}^2 \widetilde u_1(z) +2N^2\partial_{z_3}\widetilde  w_2(z)&=&\displaystyle  f_1(z')- \partial_{z_1}\widetilde p(z')&\hbox{ in }\Omega,\\
\noame
\displaystyle
-\partial_{z_3}^2   \widetilde u_2(z) -2N^2\partial_{z_3}\widetilde  w_1(z)&=&\displaystyle f_2(z')-\partial_{z_2}\widetilde p(z')&\hbox{ in }\Omega,\\
\noame
\displaystyle
-R_c\partial_{z_3}^2   \widetilde w_1(z)+4N^2 \widetilde  w_1(z) +2N^2\partial_{z_3}  \widetilde u_2(z)&=&g_1(z')\displaystyle &\hbox{ in }\Omega,\\
\noame
\displaystyle
-R_c\partial_{z_3}^2   \widetilde w_2(z)+4N^2\widetilde  w_2(z) - 2N^2\partial_{z_3} \widetilde u_1(z)&=&g_2(z')\displaystyle &\hbox{ in }\Omega,
\end{array}\right.
$$
with no-slip and no-spin condition on $\Gamma_1$ and Navier partial slip condition with no-spin condition on $\Gamma_0$.  So, we divide the proof in two steps. 

{\it Step 1.} First, we will solve the following problem
\begin{equation}\label{solving1}
\left\{\begin{array}{rcll}
\displaystyle
-\partial_{z_3}^2 \widetilde u_1(z) +2N^2\partial_{z_3}\widetilde  w_2(z)&=&\displaystyle  f_1(z')- \partial_{z_1}\widetilde p(z')&\hbox{ in }\Omega,\\
\noame
\displaystyle
-R_c\partial_{z_3}^2   \widetilde w_2(z)+4N^2\widetilde  w_2(z) - 2N^2\partial_{z_3} \widetilde u_1(z)&=&g_2(z')\displaystyle &\hbox{ in }\Omega,\\
\noame\displaystyle
\widetilde u_1(z',h(z'))=0,\quad \widetilde w_2(z',h(z'))=0,\\
\noame
\displaystyle -\partial_{z_3}\widetilde u_1(z',0)=-\lambda\widetilde u_1(z',0),\quad \widetilde w_2(z',0)=0.
\end{array}\right.
\end{equation}
From (\ref{solving1})$_1$, we have
 \begin{equation}\label{solving2}
 \partial_{z_3}\widetilde u_1(z)=(\partial_{z_1}\widetilde p(z')-f_1(z'))z_3+2N^2\widetilde w_2(z)+C(z'),
 \end{equation}
 where $C(z')$ is an unknown function. Putting this into (\ref{solving1})$_2$, we get
  \begin{equation}\label{solving3}
 \partial_{z_3}^2\widetilde w_2(z)-{4N^2\over R_c}(1-N^2)\widetilde w_2(z)=-{2N^2\over R_c}(\partial_{z_1}\widetilde p(z')-f_1(z'))z_3-{1\over R_c}g_2(z')-{2N^2\over R_c}C(z').
 \end{equation}
 The solution is
  \begin{equation}\label{solving4}
  \begin{array}{rl}
  \displaystyle
  \widetilde w_2(z)=&\displaystyle A(z')\cosh(kz_3)+B(z')\sinh(kz_3)+{1\over 2(1-N^2)}(\partial_{z_1}\widetilde p(z')-f_1(z'))z_3\\
  \noame
  &\displaystyle+{1\over 2(1-N^2)}C(z')+{1\over 4N^2(1-N^2)} g_2(z'),
  \end{array}
  \end{equation}
where $k=2N\sqrt{1-N^2\over R_c}$ and $A(z')$ and $B(z')$ are unknowns functions.

Putting this expression into (\ref{solving2}), we deduce
  \begin{equation}\label{solving5}
  \begin{array}{rl}
  \partial_{z_3}\widetilde u_1(z)=&\displaystyle {1\over 1-N^2}(\partial_{z_1}\widetilde p(z')-f_1(z'))z_3+2N^2(A(z')\cosh(kz_3)+B(z')\sinh(kz_3))\\
  \noame
  &\displaystyle
   +{1\over 1-N^2}C(z')+{1\over 2(1-N^2)}g_2(z'),
  \end{array}
  \end{equation}
  and so,
  \begin{equation}\label{solving6}
  \begin{array}{rl}
  \widetilde u_1(z)=&\displaystyle {z_3^2\over 2(1-N^2)}(\partial_{z_1}\widetilde p(z')-f_1(z'))+{2N^2\over k}(A(z')\sinh(kz_3)+B(z')\cosh(kz_3))\\
  \noame
  &\displaystyle
   +{1\over 1-N^2}C(z')z_3 +{1\over 2(1-N^2)}g_2(z')z_3+D(z').
  \end{array}
  \end{equation}
  Using the boundary conditions $\widetilde w_2(z',0)=0$ and $-\partial_{z_3}\widetilde u_1(z',0)=-\lambda\widetilde u_1(z',0)$,  we deduce
  \begin{equation}\label{CD}
  C=-2(1-N^2)A(z')-{1\over 2N^2}g_2(z'),\quad  D={2\over \lambda}(1-N^2)A(z')-{2N^2\over k}B(z')+{1\over 2\lambda N^2}g_2(z').
  \end{equation}

Then, we have
\begin{equation}\label{solving61}
  \begin{array}{rl}
  \widetilde u_1(z)=&\displaystyle {z_3^2\over 2(1-N^2)}(\partial_{z_1}\widetilde p(z')-f_1(z'))+A(z')\left({2N^2\over k}\sinh(kz_3)-2z_3+{2\over \lambda}(1-N^2)\right)\\
  \noame
  &\displaystyle +{2N^2\over k}B(z')\left(\cosh(kz_3)-1\right)-{1\over 2N^2}\left(z_3-{1\over \lambda}\right)g_2,  \\
  \\
  \displaystyle
  \widetilde w_2(z)=&\displaystyle {1\over 2(1-N^2)}(\partial_{z_1}\widetilde p(z')-f_1(z'))z_3+ A(z')(\cosh(kz_3)-1)+B(z')\sinh(kz_3).
  \end{array}
  \end{equation}
  Using the boundary condition $\widetilde w_2(z',h(z'))=0$, we deduce
\begin{equation}\label{systemQ}
\displaystyle Q_\lambda\left(\begin{array}{c} A\\
\noame
B\end{array}\right)=-{h(z')\over 2(1-N^2)}(\partial_{z_1}\widetilde p(z')-f_1(z'))\left(\begin{array}{c} h(z')\\
\noame
1\end{array}\right)+{1\over 2N^2}\left(h(z')-{1\over \lambda}\right)g_2(z')\left(\begin{array}{c} 1\\
\noame
0\end{array}\right),
\end{equation}
with
\begin{equation}\label{matQ}
Q_\lambda=\left(\begin{array}{cc} {2N^2\over k}\sinh(kh(z'))-2h(z')+{2\over \lambda}(1-N^2)& {2N^2\over k}\left(\cosh(kh(z'))-1\right)\\
\noame
\cosh(kh(z'))-1 & \sinh(kh(z'))\end{array}\right).
\end{equation}
By linearity, the solution of this system is given by
$$\begin{array}{l}
\displaystyle
A(z')=-{h(z')\over 2(1-N^2)}(\partial_{z_1}\widetilde p(z')-f_1(z'))A^\lambda_1(z')+{1\over 2N^2}\left(h(z')-{1\over \lambda}\right)g_2(z')A^\lambda_2(z'),\\
\noame
\displaystyle
B(z')=-{h(z')\over 2(1-N^2)}(\partial_{z_1}\widetilde p(z')-f_1(z'))B^\lambda_1(z')+{1\over 2N^2}\left(h(z')-{1\over \lambda}\right)g_2(z')B^\lambda_2(z'),
\end{array}$$
where $A^\lambda_1(z'), B^\lambda_1(z')$ and $A^\lambda_2(z'), B^\lambda_2(z')$ are solution of
$$Q_\lambda\left(\begin{array}{c} A^\lambda_1\\
\noame
B^\lambda_1\end{array}\right)=\left(\begin{array}{c} h(z')\\
\noame
1\end{array}\right)\quad\hbox{and}\quad Q\left(\begin{array}{c} A^\lambda_2\\
\noame
B^\lambda_2\end{array}\right)=\left(\begin{array}{c} 1\\
\noame
0\end{array}\right).$$
To compute  $A^\lambda_i, B^\lambda_i$ for $i=1,2$, we take into account that
$$\begin{array}{rl}
|Q_\lambda|=&\displaystyle{4N^2\over k}(\cosh(kh(z'))-1)-2h(z')\sinh(kh(z'))+{2\sinh(kh(z'))\over \lambda}(1-N^2)\\
\noame
=&\displaystyle 2\sinh(kh(z'))\left({2N^2\over k}{\cosh(kh(z'))-1\over \sinh(kh(z'))}-h(z')+{1\over \lambda}(1-N^2)\right),
\end{array}$$
that $Q^{-1}_\lambda$ is given by
$$
Q^{-1}_\lambda={1\over |Q_\lambda|}\left(\begin{array}{cc} \sinh(kh(z'))& -{2N^2\over k}\left(\cosh(kz_3)-1\right)\\
\noame
1-\cosh(kh(z')) & {2N^2\over k}\sinh(kh(z'))-2h(z')+{2\over \lambda}(1-N^2)\end{array}\right),
$$
and that
$$\tanh\left({kh(z')\over 2}\right)={\sinh(kh(z'))\over 1+\cosh(kh(z'))}={\cosh(kh(z'))-1\over \sinh(kh(z'))}.$$
Thus, we get that $A^\lambda_i, B^\lambda_i$, $i=1,2$, are given by (\ref{AB_partial_expression}). Then, $\widetilde u_1$ and $\widetilde w_2$ are given in (\ref{partial_u_expression}).\\

{\it Step 2.} To derive the expressions of $\widetilde u_2$ and $\widetilde w_1$, we note that $(\widetilde u_1, \widetilde w_2)$ with external forces $(f_1, g_2)$, and $(\widetilde u_2, -\widetilde w_1)$ with external forces $(f_2, -g_1)$ satisfy the same equations and boundary conditions. Then, by using this fact and the expressions of $(\widetilde u_1,\widetilde w_2)$ given in (\ref{partial_u_expression}), then we easily deduce expressions for $(\widetilde u_2,\widetilde w_1)$ given in (\ref{partial_u_expression}).

 \end{proof}

\begin{lemma}[No-slip case]\label{expression_uw_noslip}
The solution $({\bf \widetilde u}', {\bf \widetilde w}')$ of (\ref{hom_system_sub_u}) with boundary conditions (\ref{bc_limit1}) on $\Gamma_1$ and (\ref{gmenor1}) on $\Gamma_0$  are given by the following expressions:
\begin{equation}\label{perfect_u_expression}
  \begin{array}{l}
  {\bf \widetilde u}'(z)\\
  \noame ={1\over 2(1-N^2)}(\nabla_{z'}\widetilde p(z')-{\bf f}'(z'))\left[z_3^2-{h(z')}\left({2N^2\over k}\sinh(kz_3)-2z_3\right)A^\infty_1(z')-{2N^2\over k} {h(z')}\left(\cosh(kz_3)-1\right)B^\infty_1(z')
  \right]\\
  \noame
 +{1\over 2N^2} ({\bf g}(z'))^\perp\left[ z_3 - h(z')   \left({2N^2\over k}\sinh(kz_3)-2z_3\right)A^\infty_2(z')-{2N^2\over k} h(z') \left(\cosh(kz_3)-1\right) B^\infty_2(z')\right],
  \end{array}
 \end{equation}
 $$\begin{array}{l}
 {\bf \widetilde w}'(z)={1\over 2(1-N^2)}(\nabla_{z'}\widetilde p(z')-{\bf f}'(z'))^\perp\Big[z_3-h(z') (\cosh(kz_3)-1)A^\infty_1(z')-h(z')\sinh(kz_3)B_1^\infty(z')
  \Big]\\
  \noame
  +{h(z')\over 2N^2}{\bf g}'(z')\Big[  (\cosh(kz_3)-1)A^\infty_2(z') +  \sinh(kz_3)B^\infty_2(z')
  \Big],
 \end{array}
$$
where $k$ is given by (\ref{parameterk})  and  $A^\infty_i, B^\infty_2$, $i=1,2$, are given by
\begin{equation}\label{AB_partial_expression_infty}\begin{array}{ll}
\noame
A^\infty_1(z')=\displaystyle -{1\over 2},&
B^\infty_1(z')=\displaystyle {1\over 2}\coth\left({kh(z')\over 2}\right),\\
\noame
A^\infty_2(z')=\displaystyle{1\over 2}{1\over  {2N^2\over k}\tanh\left({kh(z')\over 2}\right)-h(z')},&
B^\infty_2(z')=-{1\over 2}{\tanh\left({kh(z')\over 2}\right) \over  {2N^2\over k}\tanh\left({kh(z')\over 2}\right)-h(z')}.
\end{array}
\end{equation}
\end{lemma}
\begin{proof}
The proof follows standard arguments used in the case with no-slip condition for velocity and no-spin condition for microrotation, see \cite{BayadaLuc, Luka} (see also \cite[Lemma 3.3]{Anguiano_SG_magneto} for the complete computations).

\end{proof}

\begin{remark}
We note that the expressions for $({\bf \widetilde u}', {\bf \widetilde w}')$ in the no-slip case are recovered from the expressions in the partial-slip case   considering $\lambda=+\infty$.
\end{remark}

\begin{lemma}[Perfect slip case]\label{expression_uw_noslip}
The solution $({\bf \widetilde u}', {\bf \widetilde w}')$ of (\ref{hom_system_sub_u}) with boundary conditions (\ref{bc_limit1}) on $\Gamma_1$ and (\ref{gmayor1}) on $\Gamma_0$  are given by the following expressions:
\begin{equation}\label{no_slip_u_expression}
  \begin{array}{l}
  {\bf \widetilde u'}(z) ={1\over 2(1-N^2)}(\nabla_{z'}\widetilde p(z')-{\bf f}'(z'))\left[z_3^2-h(z')^2-h(z') {2N^2\over k} \left({\cosh(kz_3)\over \sinh(h(z'))}-\coth(kh(z'))\right)\right]\\
  \noame
 -{({\bf g}'(z'))^\perp\over 2(1-N^2)}\left[z_3-h(z')+{1\over k}\left(\sinh(kz_3)-\sinh(kh(z'))\right)
+{\coth\left({kh(z')\over 2}\right)\over k}\left(\sinh(kh(z'))\cosh(kh(z'))-\cosh(kz_3)\right)
 \right],
 \end{array}
 \end{equation}
 $$\begin{array}{l}
{\bf \widetilde w}'(z)= {1\over 2(1-N^2)}(\nabla_{z'}\widetilde p(z')-{\bf f}'(z'))^\perp\left[z_3-h(z'){\sinh(kz_3)\over \sinh(kh(z'))}
\right]
+{{\bf g}'(z')\over 4N^2(1-N^2)}\left[\coth\left({kz_3\over 2}\right)-\coth\left({kh(z')\over 2}\right)\right]\sinh(kz_3),
 \end{array}
$$
%
where $k$ is given by (\ref{parameterk}).
\end{lemma}
\begin{proof}  The proof is similar to the one given in Lemma \ref{expression_uw_hat_alpha}, just taking perfect slip condition on $\Gamma_0$ instead of the partial slip one. So, we only give some details. Thus, taking into account expressions for $\widetilde w_2$ and  $\partial_{z_3}\widetilde u_1$ given in (\ref{solving4}) and (\ref{solving5}), by using the boundary conditions $w_2(z',0)=0$ and $\partial_{z_3}\widetilde u_1(z',0)=0$, we deduce
$$
A={1\over 4N^2(1-N^2)}g_2(z'),\quad C=0,
$$
which implies
\begin{equation}\label{solving6_perfecet}
  \begin{array}{rl}
  \widetilde u_1(z)=&\displaystyle {z_3^2\over 2(1-N^2)}(\partial_{z_1}\widetilde p(z')-f_1(z'))+{2N^2\over k}B(z')\cosh(kz_3) \\
  \noame
  &\displaystyle
    +{1\over 2(1-N^2)}g_2(z')\left(z_3+{1\over k}\sinh(kz_3)\right)+D(z'),\\
    \noame\displaystyle
  \widetilde w_2(z)=&\displaystyle B(z')\sinh(kz_3)+{1\over 2(1-N^2)}(\partial_{z_1}\widetilde p(z')-f_1(z'))z_3\\
  \noame
  &\displaystyle +{1\over 4N^2(1-N^2)} g_2(z')\left(\cosh(kz_3)+1\right).
  \end{array}
  \end{equation}
  Taking into account the boundary conditions $\widetilde u_1(z',h(z'))=\widetilde w_2(z',h(z'))=0$, we get
  $$Q_0\left(\begin{array}{c}
  B\\
  \noame
  D\end{array}\right)=-{h(z')\over 2(1-N^2)}(\partial_{z_1}\widetilde p(z')-f_1(z'))\left(\begin{array}{c}
  h(z')\\
  \noame
  1\end{array}\right)-{1\over 2(1-N^2)}g_2(z')\left(\begin{array}{c}
h(z')+{1\over k}\sinh(kh(z'))\\
\noame
{1\over 2N^2}\left(\cosh(kz_3)+1\right)\end{array}\right),
  $$
  with
  $$Q_0=\left(\begin{array}{cc}
{2N^2\over k}\cosh(kh(z')& 1\\
\noame
\sinh(kh(z'))&0\end{array}\right).
$$
By linearity, the solution of this system is given by
$$\begin{array}{l}
\displaystyle
B(z')=-{h(z')\over 2(1-N^2)}(\partial_{z_1}\widetilde p(z')-f_1(z'))B^0_1(z')-{1\over 2(1-N^2)}g_2(z')B^0_2(z'),\\
\noame
\displaystyle
D(z')=-{h(z')\over 2(1-N^2)}(\partial_{z_1}\widetilde p(z')-f_1(z'))D^0_1(z')-{1\over 2(1-N^2)}g_2(z')D^0_2(z'),
\end{array}$$
where $B^0_1(z'), D^0_1(z')$ and $B^0_2(z'), D^0_2(z')$ are solution of
$$Q_0\left(\begin{array}{c} B^0_1\\
\noame
D^0_1\end{array}\right)=\left(\begin{array}{c} h(z')\\
\noame
1\end{array}\right)\quad\hbox{and}\quad Q_0\left(\begin{array}{c} B^0_2\\
\noame
D^0_2\end{array}\right)=\left(\begin{array}{c}
h(z')+{1\over k}\sinh(kh(z'))\\
\noame
{1\over 2N^2}\left(\cosh(kh(z'))+1\right)\end{array}\right).$$
To compute  $B^0_i, D^0_i$ for $i=1,2$, we take into account that
$$\begin{array}{rl}
|Q_0|=&\displaystyle-\sinh(kh(z')),
\end{array}$$
and that $Q^{-1}_0$ is given by
 $$Q_0^{-1}=\left(\begin{array}{cc}
0& {1\over \sinh(kh(z'))}\\
\noame
1&-{2N^2\over k}\cosh(kh(z'))\end{array}\right).
$$
Thus, we get that $A^\lambda_i, B^\lambda_i$, $i=1,2$, are given by
$$\begin{array}{ll}
\displaystyle B_1^0(z')={1\over \sinh(kh(z'))},&\displaystyle D_1^0(z')=h(z')-{2N^2\over k}\cosh(kh(z')),\\
\noame
\displaystyle
B_2^0(z')={1\over 2N^2}\coth\left({kh(z')\over 2}\right),&\displaystyle D_2^0(z')=h(z')+{1\over k}\sin(kh(z'))-{1\over k}\left(1+\cosh(kh(z')))\cosh\left(kh(z')\right)\right).
\end{array}$$
Then, $\widetilde u_1$ and $\widetilde w_2$ are given by (\ref{no_slip_u_expression}). As in the case of partial slip condition, Then,  $\widetilde u_2$ and $\widetilde w_1$ are deduced from expressions of $\widetilde u_1$ and $\widetilde w_2$, and are given  in (\ref{no_slip_u_expression}). 
\end{proof}

\section{Reynolds equation}\label{sec:Reynolds}
 In this section we given the main result of the paper, which is the generalized 2D Reynolds equations for the pressure in every case.  To do this, we proceed as follows. First, we compute the average limit velocity 
 $${\bf \widetilde U}(z')=\int_0^{h(z')}{\bf \widetilde u}(z)\,dz_3,$$
 with the expressions given for $\widetilde u$, and then, we use the divergence condition (\ref{divzh_limit}) to derive the Reynolds equation for pressure.

 \begin{theorem}[Reynolds equation in the partial slip case]The limit pressure $\widetilde p\in L^2_0(\Omega)\cap H^1(\omega)$ verifies the following variational equation:
\begin{equation}\label{Reynolds_partial}\begin{array}{l}
\displaystyle\int_\omega{h^3(z')\over 1-N^2}\Theta^\lambda_1(z')\nabla_{z'}\widetilde p(z')\cdot \nabla_{z'}\theta(z')\,dz'-{1\over \lambda}\int_\omega h^2(z')\Theta^\lambda_3(z')\nabla_{z'}\widetilde p(z')\cdot \nabla_{z'}\theta(z')\,dz'\\
\noame
\displaystyle
=\int_\omega{h^3(z')\over 1-N^2}\Theta^\lambda_1(z'){\bf f}'(z')\cdot \nabla_{z'}\theta(z')\,dz'-\int_\omega {h(z')\over 2N^2}\Theta^\lambda_2(z')({\bf g}'(z'))^\perp\cdot \nabla_{z'}\theta(z')\,dz'\\
\noame
\displaystyle+{1\over \lambda}\left(-\int_\omega {h^2(z')}\Theta_3^\lambda(z'){\bf f}'(z')\nabla_{z'}\theta(z')\,dz'+{1\over 2N^2}\int_\omega\Theta^\lambda_4(z')({\bf g}'(z'))^\perp\cdot \nabla_{z'}\theta(z')\,dz'
\right)\\
\noame
\displaystyle
-{1\over \lambda^2}\int_\omega{h(z')(1-N^2)\over N^2}\Theta^\lambda_5(z')({\bf g}'(z'))^\perp\cdot \nabla_{z'}\theta(z')\,dz',
\end{array}
\end{equation}
for any $\theta\in H^1(\omega)^2$, where
$$\begin{array}{l}
\Theta_1^\lambda(z')=
 {1\over 6}+{1\over 2}A_1^\lambda(z')+{N^2\over kh(z')}B_1^\lambda(z')-{N^2\over k^2h^2(z')}\left((\cosh(kh(z'))-1)A_1^\lambda(z')+\sinh(kh(z'))B_1^\lambda\right),\\
\noame
\Theta_2^\lambda(z')=   {h(z')\over 2}-  \left({2N^2\over k}(\cosh(kh(z'))-1)-h(z')^2\right)A^\lambda_2(z') -{2N^2\over k}   \left({1\over k}\sinh(kh(z'))-h(z')\right) B^\lambda_2(z'),\\
\noame
\Theta_3^\lambda(z')= A^\lambda_1(z'),\\
 \noame
\Theta_4^\lambda(z')=
 1-\left(-2h(z')(1-N^2)+{2N^2\over k^2}(\cosh(kh(z'))-1)-h(z')^2\right)A^\lambda_2(z') -{2N^2\over k}\left({1\over k}\sinh(kh(z'))-h(z')\right)B^\lambda_2(z'),
 \\
 \noame
 \Theta_5^\lambda(z')=A^\lambda_2(z'),
\end{array}$$
with $k$ given in (\ref{parameterk}) and $A_i^\lambda, B_i^\lambda$, $i=1,2$, given in (\ref{AB_partial_expression}).
%
%
 \end{theorem}
\begin{proof}Integrating with respect to $z_3$ between $0$ and $h(z')$ the expressions of ${\bf \widetilde u}'$ given in (\ref{partial_u_expression}):
$$
\begin{array}{rl}
{\bf \widetilde U}'(z')
= &\displaystyle {h^3(z')\over 1-N^2}\Theta_1^\lambda(z')(\nabla_{z'}\widetilde p(z')-{\bf f}'(z')) +{h(z')\over 2N^2}\Theta_2^\lambda (z')({\bf g}'(z'))^\perp\\
\noame
&\displaystyle +{1\over \lambda}\left[-h^2(z')\Theta_3^\lambda (z')(\nabla_{z'}\widetilde p(z')-{\bf f}'(z'))-{1\over 2N^2}\Theta_4^\lambda(z')({\bf g}'(z'))^\perp\right]+{1\over \lambda^2}{h(z')(1-N^2)\over N^2}\Theta_5^\lambda(z')({\bf g}'(z'))^\perp,
%
%
\end{array}$$
where $\Theta_i^\lambda$, $i=1, \ldots, 5$ are given in the theorem.

Taking these expressions into the divergence condition (\ref{divzh_limit}), we derive the Reynolds equation in the variational form (\ref{Reynolds_partial}).

\end{proof}

 \begin{theorem}[Reynolds equation in the no-slip case]The limit pressure $\widetilde p\in L^2_0(\Omega)\cap H^1(\omega)$ verifies the following equation:
\begin{equation}\label{Reynolds_noslip}\begin{array}{l}
\displaystyle\int_\omega{h^3(z')\over 1-N^2}\Phi(z')\nabla_{z'}\widetilde p(z')\cdot \nabla_{z'}\theta(z')\,dz'
=\int_\omega{h^3(z')\over 1-N^2}\Phi(z'){\bf f}'(z')\cdot \nabla_{z'}\theta(z')\,dz'
\end{array}
\end{equation}
for any $\theta\in H^1(\omega)^2$, where
\begin{equation}\label{Phifunc}\Phi(h(z'),N,R_c)={1\over 12}+{R_c\over 4h(z')^2(1-N^2)}-{1\over 4h(y')}\sqrt{N^2R_c\over 1-N^2}\coth\left(Nh(z')\sqrt{1-N^2\over R_c}\right).
\end{equation}
\end{theorem}
 \begin{proof}
This is the classic setting, so we give a some details. Integrating  with respect to $z_3$ between $0$ and $h(z')$ the expression of ${\bf \widetilde u}'$ given in  (\ref{no_slip_u_expression}):
$$\begin{array}{rl}{\bf \widetilde U}'(z')
  =&-{h^3(z')\over 1-N^2}\Phi(z')(\nabla_{z'}\widetilde p(z')-{\bf f}'(z'))+{1\over N^2}\Psi(z')({\bf g}'(z'))^\perp,
  \end{array}
$$
where, taking into account the definition of $k$ given in (\ref{parameterk}), we get that  $\Phi$ is given by (\ref{Phifunc}) and
$$\begin{array}{rl}\Psi(z')= &\displaystyle {h(z')^2\over 2}- {h(z')\over 2}   \left({2N^2\over k^2}(\cosh(kh(z'))-1)-h(z')^2\right){1\over {2N^2\over k}\tanh\left({kh(z')\over 2}\right)-h(z')}\\
\noame
&\displaystyle +{N^2\over k} h(z') \left({1\over k}\sinh(kh(z'))-h(z')\right) {\tanh\left({kh(z')\over 2}\right)\over {2N^2\over k}\tanh\left({kh(z')\over 2}\right)-h(z')}=0.
\end{array}$$
 \end{proof}

  \begin{theorem}[Reynolds equation in the perfect slip case]The limit pressure $\widetilde p\in L^2_0(\Omega)\cap H^1(\omega)$ verifies the following variational equation:
\begin{equation}\label{Reynolds_perfect}\begin{array}{rl}
\displaystyle \int_\omega{h^3(z')\over 3(1-N^2)}\Theta^0_1(z')\nabla_{z'}\widetilde p(z')\cdot \nabla_{z'}\theta(z')\,dz'
=&\displaystyle \int_\omega{h^3(z')\over 3(1-N^2)}\Theta^0_1(z'){\bf f}'(z')\cdot \nabla_{z'}\theta(z')\,dz'\\
\noame&\displaystyle+\int_\omega{h^2(z')\over 4(1-N^2)}\Theta^0_2(z')({\bf g}'(z'))^\perp\cdot \nabla_{z'}\theta(z')\,dz',
\end{array}
\end{equation}
for any $\theta\in H^1(\omega)^2$, where
$$\begin{array}{l}
\Theta_1^0(z')= 1+ {3N^2\over k h^2(z')} \left({1\over k}-\coth(kh(z'))\right),\\
\noame
\Theta_2^0(z')= 1-
 {2\over kh^2(z')}\sinh(kh(z'))\left(-1-{1\over k}{1\over \cosh(kh(z'))\sinh(kh(z'))}+\cosh(kh(z'))\coth\left({kh(z')\over 2}\right)\right),
\end{array}$$
with $k$ given by (\ref{parameterk}).
%
%
 \end{theorem}
 \begin{proof}Integrating with respect to $z_3$ between $0$ and $h(z')$ the expressions of ${\bf \widetilde u}'$ given in (\ref{perfect_u_expression}):
$$
\begin{array}{l}
\displaystyle {\bf \widetilde U}'(z')
 = -{h^3(z')\over 3(1-N^2)}\Theta_1^0(z')(\nabla_{z'}\widetilde p(z')-{\bf f}'(z'))+{h^2(z')\over 4(1-N^2)}\Theta_2^0(z')({\bf g'}(z'))^\perp,
 \end{array}$$
where $\Theta_i^0$, $i=1,2,$ are given in the theorem.

Tanking these expressions into the divergence condition (\ref{divzh_limit}), we derive the Reynolds equation in the variational form (\ref{Reynolds_perfect}).

\end{proof}
\section*{Acknowledgements}
The first author  belongs to the ``Mathematical Analysis" Research Group (FQM104) at Universidad de Sevilla.
The second author is supported by the Croatian Science Foundation under the project AsyAn (IP-2022-10-1091) and Republic of Croatia's MSEY in course of Multilateral scientific and technological cooperation in Danube region under the project MultiHeFlo.

\section*{Declarations}
\paragraph{Conflict of interest} The authors declared that they have no conflict of interest.

\paragraph{Data availability statement} Data sharing not applicable - no new data generated.

\end{document}